\documentclass[11pt,reqno]{amsart}
\usepackage[left=33mm,right=33mm,top=30mm,bottom=32mm]{geometry}
\usepackage{mathtools,amssymb,amsthm,mathrsfs,color,float}

\usepackage[colorlinks,
linkcolor=red,
anchorcolor=green,
citecolor=blue, 
]{hyperref}

\usepackage[T1]{fontenc}
\usepackage[utf8]{inputenc}

\usepackage{calc}
\definecolor{bleu1}{RGB}{0,57,128}
\def\bleu1{\color{bleu1}}

\usepackage{etoolbox}
\patchcmd{\section}{\normalfont}{\normalfont \bleu1}{}{}
\patchcmd{\subsection}{\normalfont}{\normalfont \bleu1}{}{}
\patchcmd{\subsubsection}{\normalfont}{\normalfont \bleu1}{}{}

\newtheorem{THEalpha}{\bleu1 Theorem}%[section]
\newtheorem{CORalpha}[THEalpha]{\bleu1 Corollary} 

\newtheorem{The}{\bleu1 Theorem}[section]
\newtheorem{Pro}{\bleu1 Proposition}[section]
\newtheorem{Lem}{\bleu1 Lemma}[section]
\newtheorem{Cor}{\bleu1 Corollary}[section]

\theoremstyle{definition}
\newtheorem{defn}{\bleu1 Definition}[section]
\newtheorem*{Prbm}{Problem}
\newtheorem{Rem}{\bleu1 Remark}[section]
\newtheorem{Ex}{\bleu1 Example}[section]

\newcommand{\T}{\mathbb{T}}
\newcommand{\R}{\mathbb{R}}
\newcommand{\Z}{\mathbb{Z}}

\newcommand{\LL}{\mathbf{L}}
\newcommand{\HH}{\mathbf{H}}
\newcommand{\lam}{\lambda}
\newcommand{\supp}{\operatorname{supp}}
\newcommand{\wtM}{\widetilde{\mathfrak{M}}}
\newcommand{\tmu}{\tilde\mu}

\newcommand{\tmulax}{\tilde\mu_x^\lambda}
\newcommand{\xilax}{\xi_x^\lambda}
\newcommand{\dxilax}{\dot{\xi}_x^\lambda}

\newcommand{\Psma}{\mathbb{P}^\sigma_{_{G}}}

\newcommand{\Psmone}{\mathbb{P}_{_{G}}}

\newcommand{\mfK}{\mathfrak{K}}
\newcommand{\mfD}{\mathfrak{D}}

\def\leq{\leqslant}
\def\geq{\geqslant}
\def\tilde{\widetilde}
\def\du#1{\langle#1\rangle}
\def\cA{\mathcal{A}}

\title[]{The selection problem for a new class of perturbations of Hamilton-Jacobi equations and its applications}

\author{Qinbo Chen}
\address[Qinbo Chen]{School of Mathematics, Nanjing University, Nanjing 210093, China}
\email{qinbochen@nju.edu.cn}

\begin{document}

\begin{abstract}
This paper studies a perturbation problem given by the equation:
\begin{equation*}
H(x, d_xu_\lambda, \lambda u_\lam(x))+\lambda V(x,\lambda)=c  \quad \text{in $M$},
\end{equation*}
where $M$ is a closed manifold and $\lambda>0$ is a perturbation parameter. The Hamiltonian $H(x,p,u):T^*M\times\R\to \R$ satisfies certain convexity, superlinearity, and monotonicity conditions.  $\lambda V(\cdot,\lambda):M\to\R$ converges to zero as $\lambda\to 0$. First, we study the asymptotic behavior of the viscosity solution $u_\lambda:M\to\R$ as $\lambda$ approaches zero. This perturbation problem explores the combined effects of both the vanishing discount process and potential perturbations, leading to a new selection principle that extends beyond the classical vanishing discount approach. Additionally, we apply this principle to Hamilton-Jacobi equations with $u$-independent Hamiltonians, resulting in the introduction of a new solution operator. This operator provides new insights into the variational characterization of viscosity solutions and Mather measures.

\end{abstract}

\maketitle

%\tableofcontents

\section{Introduction}\label{section_Int}

This paper addresses a class of perturbation problems within the framework of first- order Hamilton-Jacobi equations. Consider the perturbation problem defined by the equation
\begin{equation}\label{introeq1}
H(x, d_xu_\lambda, \lambda u_\lam(x))+\lambda V(x,\lambda)=c  \quad \text{in $M$},
\end{equation}
where $\lambda>0$ serves as a perturbation parameter and $M$ is a closed, connected Riemannian manifold. The function $H(x,p,u):T^*M\times\R\to \R$ is a given continuous Hamiltonian, convex and superlinear in  the momentum variable $p$, and is strictly increasing in $u$. A typical example is the discounted system $H(x,p,u)=G(x,p)+u$. The potential $V(\cdot,\lambda)$ is a family of continuous functions on $M$, with $\lambda V(\cdot,\lambda)$ converging to zero as $\lambda\to 0$. In this setting, equation \eqref{introeq1} is an approximation 
    of the critical Hamilton-Jacobi equation: 
\begin{equation}\label{introeq2}
H(x, d_xu, 0)=c  \quad \text{in $M$},
\end{equation}
where $c$ is the \textit{critical value} or \textit{ergodic constant}. This value $c$ is the unique constant for which the equation \eqref{introeq2} has solutions. The notion of solution here is understood in the viscosity sense as introduced by Crandall and Lions \cite{Crandall_Lions1983}.

A key question is to \textit{understand the asymptotic behavior of the viscosity solution $u_\lambda:M\to\R$ of the perturbation problem \eqref{introeq1} as $\lambda$ tends to zero}.  In the special case where the potential term is absent (i.e. $V\equiv 0$), this problem reduces to the selection problem for discounted approximations, known as the {\em vanishing discount} problem, which has been extensively studied due to its connections to fields such as homogenization of HJEs and optimal control theory, see for instance \cite{Lions_Papanicolaou_Varadhan1987, Bardi_Capuzzo-Dolcetta1997, Arisawa_Lions1998, DFIZ2016, Tran_book2021} and the references therein. We also refer to recent works \cite{Arnaud_Viterbo2024higher, Arnaud_Su, Maro_Sorrentino2017} for dynamical and geometric perspectives on the vanishing discount problem.

However, in the non-zero potential case (i.e., $V\not\equiv 0$), the corresponding selection problem seems more complex, as it involves the {\em combined effects of both the vanishing discount process and potential perturbations}. This new approximation scheme is interesting and useful, as {\em the vanishing discount method alone cannot select all solutions of the critical Hamilton-Jacobi equation \eqref{introeq2}.} This limitation highlights the necessity for alternative selection principles, motivating our investigation into the new perturbation problem \eqref{introeq1}.

Let us briefly review the existing literature on related problems. In the case where the potential is zero ($V\equiv 0$), which corresponds to the vanishing discount problem, the study of the asymptotic convergence of viscosity solutions  traces back to the seminal work of Lions, Papanicolaou, and Varadhan \cite{Lions_Papanicolaou_Varadhan1987}.  They introduced the discount approximation:
\begin{equation}\label{introeq_4}
	\lambda w_\lambda+G(x, d_xw_\lambda)=0,
\end{equation}
for coercive Hamiltonians $G$, and showed that $\lambda w_\lambda$ converges to  $-c(G)$, where $c(G)$ is the critical value. Also, the sequence $\{w_\lambda+c(G)/\lambda\}_{\lambda>0}$ has subsequences converging to viscosity solutions of the corresponding ergodic Hamilton-Jacobi equation
\begin{equation}\label{introeq_5}
	G(x, d_xw)=c(G).
\end{equation}
However, due to the non-uniqueness of solutions to equation \eqref{introeq_5}, it was unclear at the time whether the limits of the sequence along different subsequences would yield the same solution to \eqref{introeq_5}. Partial progress was made in \cite{Gomes2008, Iturriaga_Sanchez-Morgado2011}, but a complete convergence result was first established by \cite{DFIZ2016} utilizing tools from weak KAM theory. Since then, the selection problem has been successfully dealt with under various settings and through different approaches. These include generalizations to second order settings \cite{Mitake_Tran2017} (degenerate viscous case) \cite{Ishii_Mitake_Tran2017_1, Ishii_Mitake_Tran2017_2} (fully nonlinear case), weakly coupled systems \cite{Davini_Siconolfi_Zavidovique_18, Davini_Zavidovique2017, Ishii2019vanishing, Ishii_Jin2020}, mean field games \cite{Cardaliaguet_Porretta_MFG19},  negative discounting \cite{Davini_Wang2020,WYZ_negative},  non-compact domains \cite{Ishii_Siconolfi2020, Davini2022}, and discounted state-constraint equations  \cite{SonTu2022}. All the results aforementioned require convexity of the Hamiltonians, some possible non-convex cases were studied in \cite{Gomes_Mitake_Tran2018}. Moreover, some recent works \cite{CCIZ2019, Gomes_Mitake_Tran2018, WYZ_2021, Zavidovique2022_degenerate, QC, Zhang2024limit, Chen_Fathi_Za_Zha2024,DNYZ2024} have also addressed selection problems for nonlinear discounted equations and contact Hamiltonians, further enriching the field.

In contrast to the well-studied zero potential case, the non-zero potential case has not yet been thoroughly investigated. This gap in the literature motivates our current study. Our work provides a new selection principle that extends beyond the vanishing discount process, offering the flexibility to select a broader variety of solutions to the critical Hamilton-Jacobi equation. Furthermore, we demonstrate how this principle provides new insights into classical theories, including weak KAM theory and Aubry-Mather theory.

More precisely, the objectives of this paper are twofold:
\begin{enumerate}
	\item [\bf(I)] We establish the convergence of the solution $u_\lambda$ to the perturbation problem \eqref{introeq1} as $\lambda\to 0$, thereby selecting a particular solution to the critical Hamilton-Jacobi equation (see Theorem \ref{mainresult1} and Theorem \ref{mainresult2}). Notably, by choosing different potential functions $V$, one can select a wider range of solutions for the critical equation, extending beyond the limitations of the vanishing discount method.

\item [\bf(II)] We explore the applications of the new selection principle in the context of classical first-order Hamilton-Jacobi equations with $u$-independent convex Hamiltonians $G(x,p):T^*M\to\R$. We demonstrate how this principle facilitates the introduction of a novel solution operator (Definition \ref{def_operaP}), which offers a new variational representation of viscosity solutions (Theorem \ref{mainresult3}). The analysis of this operator reveals that by choosing suitable potential perturbations, our approximation scheme enables the selection of \emph{any desired viscosity solution} for $G$.  Additionally, it can be used to describe the uniqueness structure of the solutions (Theorem \ref{mainresult4}), and the structural properties of equilibrium Mather measures (see Theorem \ref{mainresult6}).
\end{enumerate}

\section{Main results}\label{section_MainR} 

\subsection{Part I: Asymptotic convergence}
\subsubsection{Hypotheses}
Let $M$ be a closed connected smooth manifold. Denote by $T^*M$  the cotangent bundle of $M$, and by $(x,p)$ a point of $T^*M$. The continuous Hamiltonian $H(x,p,u):T^*M\times\R\to \R$ used in this paper satisfies the following conditions: 
\begin{itemize}
\item[\bf(H1)] (Convexity)  $H(x,p,u)$ is convex in $p$.
\item[\bf(H2)]  (Superlinearity) 
$H(x,p,u)/\|p\|_x \to +\infty \quad \text{as~ $\|p\|_x\to +\infty$.} $

\item[\bf (H3)] (Monotonicity)  $ H(x,p,u)$ is strictly increasing in $u$.
\item[\bf (H4)] ($u$-derivative at $0$)  The partial derivative $\frac{\partial H}{\partial u}(x,p,0)$ exists and  $\frac{\partial H}{\partial u}(x,p,0)>0$ for every $(x,p)\in T^*M$. Moreover, for every compact subset $S\subset T^*M$, we can find  a modulus of continuity\footnote{A modulus of continuity is a non-decreasing function $\eta_{_S}:[0,+\infty)\to[0,+\infty)$ with $\lim_{u\to 0}\eta_{_S}(u)=0$.} $\eta_{_S}$ such that
 \begin{equation}\label{reuwe2o}
 \left|H(x,p,u)-H(x,p,0)-u\frac{\partial H}{\partial u}(x,p,0)\right|\leq |u| \eta_{_S}(| u|),\quad \forall~(x,p,u)\in S\times \R.
 \end{equation}
 It is worth noting that the requirement \eqref{reuwe2o} always holds when $H$ is of class $C^1$. 
\end{itemize}

\,

A good example to keep in mind is the dissipative mechanical system $H(x,p,u)=\|p\|^2+U(x)+\sigma(x)u$, with $\sigma>0$.
 The convexity condition {\bf(H1)} and the superlinearity condition {\bf(H2)} are standard in weak KAM theory and Aubry-Mather theory. The monotonicity condition {\bf(H3)} ensures that the associated equation obeys the comparison principle.    
\subsubsection{The Lagrangian}  
The Hamiltonian $H$ has a conjugated Lagrangian $L: TM\times \R\to \R$, defined via the Fenchel formula 
 \begin{align*}
	L(x, v, u)=\sup_{p\in T^*_x M}\Big(\du{p,v}_x-H(x,p, u)\Big), \quad\text{for every $(x,v,u)\in TM\times\R$}.
\end{align*}
Due to conditions {\bf(H1)}--{\bf(H4)} on $H$, it is easy to check that  $L$ is continuous and satisfies properties analogous to those of $H$:
\begin{itemize}
\item  
 $L(x,v,u)$ is convex and superlinear in the fibers, and is strictly decreasing  in $u$.

\item 
 The partial derivative  $\frac{\partial L}{\partial u}(x,v,0)$ exists and $\frac{\partial L}{\partial u}(x,v,0)<0$ for every $(x,v)\in TM$. Moreover, for every compact subset $S\subset TM$, there is  a modulus of continuity $\eta_{_S}$ such that
 \begin{equation}\label{d23rL4}
 \left|L(x,v,u)-L(x,v,0)-u\frac{\partial L}{\partial u}(x,v,0)\right|\leq |u| \eta_{_S}(| u|),\quad \forall~(x,v,u)\in S\times \R.
 \end{equation}
\end{itemize} 

\,

For brevity, we can introduce the Lagrangian $L^r:TM\to \R, (x,v)\mapsto  L(x,v,r)$ for each fixed $r\in \R$, and its associated Hamiltonian $H^r:T^*M\to\R$, $(x,p)\mapsto H(x,p,r)$.

\subsubsection{The convergence result}
For $\lambda>0$, we consider the following perturbation problem: 
 \begin{equation}\label{eq_H_lambda}\tag{P$_\lambda$}
H(x,d_xu_\lambda, \lambda u_\lam(x))+\lambda V(x,\lambda)=c(H^0)  \quad \text{in $M$}
\end{equation} 
where the potential $ V(\cdot,\lambda): M\to \R$ is continuous, and $c(H^0)$ denotes the critical value of the Hamiltonian $H^0(x,p)=H(x,p,0)$. 
Our first main result concerns the asymptotic behavior of the viscosity solution $u_\lambda$ of \eqref{eq_H_lambda} as $\lambda\to 0^+$:

\begin{THEalpha}\label{mainresult1}
Let $H: T^*M\times\R\to \R$ satisfy conditions {\bf(H1)}--{\bf(H4)}, and suppose that $V(\cdot,\lambda)$ converges uniformly to $V(\cdot,0)$ as $\lambda\to 0$. Then there exists $\lambda_0>0$ such that for every $\lambda\in (0,\lambda_0]$, equation \eqref{eq_H_lambda} has a unique continuous viscosity solution $u_\lambda: M\to \R$. Moreover, $u_\lambda$ uniformly converges, as $\lambda\to 0$, to a particular  function $u_0: M\to \R$ which is a viscosity solution of the critical Hamilton-Jacobi equation 
\begin{equation}\label{eq_G}\tag{E}
H(x,d_xu, 0)=c(H^0)  \quad \text{in $M$}.
\end{equation}
Specifically, $u_0$ is the largest viscosity subsolution $w$ of \eqref{eq_G} such that
\begin{equation*}
\int_{TM}w(x)\frac{\partial L}{\partial u }(x,v,0) \, d\tmu(x,v)\geq  \int_{M}V(x,0) \, d\mu(x),
\end{equation*}
for any Mather measure $\tmu$ of the Lagrangian $L^0$.
\end{THEalpha}

\begin{Rem}
For the concept of Mather measures see Definition \ref{Def_Mather}. In short, Mather measures are action minimizing measures on $TM$ for the Lagrangian system.  For any Mather measure $\tilde\mu$ on $TM$, we denote by $\mu=\pi_\# \,\tilde\mu$ the corresponding projected Mather measure on $M$. 
\end{Rem}
\begin{Rem} It is important to note that:
\begin{enumerate}
\item The assumption that  $\lim_{\lambda\to 0}V(\cdot,\lambda)$ exists is necessary, otherwise, the solutions of \eqref{eq_H_lambda} may not converge, as illustrated by the following counterexample  
	\[\lambda u_\lambda(x)+\|d_xu_\lambda\|^2-\cos^2x+\lambda V(x,\lambda)=0\quad\text{in $\T$},\]
	with $V(x,\lambda)=\sin x-\lambda^{-1/2}$. In this case, $u_\lambda(x)=-\sin x+\lambda^{-1/2}$ fails to converge.
	
		\item 
 Equation \eqref{eq_H_lambda} may not admit solutions for large $\lambda>\lambda_0$, see Example \ref{ex_hasnosol}. Moreover, by using standard arguments as in \cite{DFIZ2016}, the superlinearity {\bf(H2)} can be relaxed to a coercivity condition.  However, as shown in \cite{Ziliotto_counterexample_2019, NiPanrui_Multiple}, the convergence result may not hold if either the condition {\bf{(H1)}} or the condition {\bf(H3)} is not satisfied. 
\end{enumerate}
	\end{Rem}

Our next theorem provides another formula for the limit $u_0$. It involves the Peierls barrier (see Definition \ref{dsadw14}), Mather measures, and the potential $V(x,0)$.  
 
\begin{THEalpha}\label{mainresult2}
	The limit solution $u_0: M\to \R$ obtained in Theorem  \ref{mainresult1}
	satisfies:  
	\begin{equation}\label{particular_solution}
		u_0(x)=\inf_{\tilde{\mu}\in \wtM(L^0)}  \frac{\int_{TM}  h(y,x)\frac{\partial L}{\partial u}(y,v,0)\,d\tilde\mu(y,v)+\int_M V(y,0)\,d\mu(y)}{\int_{TM} \frac{\partial L}{\partial u}(y,v,0) \,d\tilde\mu(y,v)   },
    \end{equation}
where $h: M\times M\to\R $ is the Peierls barrier of the Lagrangian $L^0$, and $\wtM(L^0)$ denotes the family of all Mather measures of the Lagrangian $L^0$.\end{THEalpha}

We remark that  for any $x\in M$, this infimum is attained since the set $\wtM(L^0)$ is compact.

\begin{Rem}[Comparison with vanishing discount limit]
Theorem \ref{mainresult2} shows that our approximation scheme \eqref{eq_H_lambda} indeed provides a new selection principle for viscosity solutions or weak KAM solutions.  Specifically, it is known (see for instance \cite{DFIZ2016, Chen_Fathi_Za_Zha2024}) that the classical vanishing discount approach selects a particular solution $v_0$ of the form
\begin{equation*}
		v_0(x)=\inf_{\tilde{\mu}\in \wtM(L^0)}  \frac{\int_{TM}  h(y,x)\frac{\partial L}{\partial u}(y,v,0)\,d\tilde\mu(y,v)}{\int_{TM} \frac{\partial L}{\partial u}(y,v,0) \,d\tilde\mu(y,v)   }.
    \end{equation*}
However, such discount approximations are limited in their ability to select all possible solutions of the critical H-J equation, regardless of how the term $\partial L/\partial u$ is adjusted. In contrast, the solution \eqref{particular_solution} selected by our approximation scheme \eqref{eq_H_lambda} offers more flexibility. By varying the potential term $V$, we are able to select a wider variety of solutions to the critical H-J equation. This flexibility emphasizes the utility of the new selection principle, expanding the range of possible solutions beyond what the vanishing discount approach offers. Actually, it enables us to select \textbf{any desired viscosity solution} to the critical equation (as indicated by  Corollary \ref{mainre31} in the subsequent section \ref{sub_applickqw}). See Example \ref{ex_simexamp2} for a simple example illustrating this phenomenon. 
	
\end{Rem}

\subsection{Part II: Application}\label{sub_applickqw}

\,

We explore the applications of the results established above to study classical first-order convex Hamilton-Jacobi equations with $u$-independent system. 
Consider a continuous Hamiltonian  $G(x,p):T^*M\to\R$ that satisfies:
\begin{enumerate}
	\item [\bf(G1)] $G(x,p)$ is convex in $p$.
	
	\item [\bf(G2)]  $G(x,p)/\|p\|_x\to +\infty$ as $\|p\|_x\to \infty$.
\end{enumerate} 
The stationary Hamilton-Jacobi equation associated to $G$ is given by:
\begin{equation}\label{eq_classical}\tag{G}
	G(x, d_x u)=c(G) \quad\text{in $M$}
\end{equation}
where $c(G)$ represents the critical value of $G$.

\subsubsection{A new characterization for solutions of \texorpdfstring{\eqref{eq_classical}}{}}
In classical weak KAM theory, it is known  (see \cite[Theorem 0.2]{Contreras_action_potential} or \cite{Fathi_book}) that any weak KAM solution to equation \eqref{eq_classical} can be characterized using the projected Aubry set and the Ma\~n\'e potential. Here, we propose a new characterization for all viscosity solutions of \eqref{eq_classical} by utilizing Mather measures and the Peierls barrier (see Theorem \ref{mainresult3} below). We need to fix some notations.

\begin{itemize}
	
	\item The Lagrangian associated with the Hamiltonian $G$ is denoted as
$L_G(x,v): TM\to\R,$ which is continuous on $TM$ and is convex and superlinear in $v$.
\item Denote by $ \tilde{\mathfrak{M}}(L_G)$  the set of all  Mather measures of $L_G$, and by $ \mathfrak{M}(L_G)$ the set of all  projected Mather measures of $L_G$.
\item $h_G(\cdot,\cdot): M\times M\to \R$ ~---~ the Peierls barrier of $L_G$.
\item $C(M)$ ~---~ the space of all continuous functions from $M$ to $\R$, and we endow, as usual, this space with the sup norm $\|f\|_\infty=\sup_{x\in M}|f(x)|$.
\end{itemize}

 Inspired by Theorem \ref{mainresult2}, we define the following:
\begin{defn}\label{def_operaP}
	Given a continuous function $\sigma: M\to (0,+\infty)$, we define an operator \[\Psma: C(M)\to  C(M)\] 
	by setting for every $\varphi\in C(M)$ the function $\Psma\varphi:M\to \R$ as follows
\begin{equation}\label{dfpo01}
\begin{aligned}
	\Psma\varphi(x):=&\inf_{\tmu \in \tilde{\mathfrak{M}}(L_G)} \left\{ \frac{\int_{TM}  \sigma(y) h_G(y,x)\,d\tmu(y,v)}{\int_{TM} \sigma(y) \,d\tmu(y,v)   } +\frac{\int_{TM}  \sigma(y)\varphi(y)\,d\tmu(y,v)}{\int_{TM} \sigma(y) \,d\tmu(y,v)} \right\}\\
	=&\inf_{\mu\in \mathfrak{M}(L_G)} \left\{ \frac{\int_{M}  \sigma(y) h_G(y,x)\,d\mu(y)}{\int_{M} \sigma(y) \,d\mu(y)   } +\frac{\int_M  \sigma(y)\varphi(y)\,d\mu(y)}{\int_{M} \sigma(y) \,d\mu(y) } \right\}.
\end{aligned}
    \end{equation}
\end{defn}

\,

A special example to keep in mind is to set $\sigma\equiv 1$ on $M$.

The next theorem below indicates that  $\Psma$ is actually a solution operator. As we will show in Section \ref{section_proofofCandD}, this relies on applying Theorem \ref{mainresult1} and Theorem \ref{mainresult2} to a specific perturbation problem given by
\begin{equation*}
	\sigma(x)\lambda u_\lambda(x)+G(x, d_xu_\lambda)-\lambda V(x,\lambda)=c(G)\quad \text{in $M$}
\end{equation*}
where $V(x,\lambda)$ converges,  as $\lambda\to 0$, to $V(x,0)=\sigma(x)\varphi(x)$.

\begin{THEalpha}\label{mainresult3}
Suppose that $G: T^*M\to \R $ satisfies {\bf(G1)}--{\bf(G2)}. Fix a continuous function $\sigma: M\to (0,+\infty)$. Then, the operator $\Psma$ has the following properties:
    \begin{itemize}
        \item [(i)] \textbf{Lipschitz continuity}: $\Psma$ is Lipschitz continuous with Lipschitz constant $1$, i.e., 
        \begin{equation*}
        	\|\Psma \varphi_1-\Psma\varphi_2 \|_\infty\leq \|\varphi_1-\varphi_2\|_\infty,\quad \text{for any $\varphi_1, \varphi_2\in C(M)$}.
        \end{equation*}
    	
    	\item [(ii)] \textbf{Solutions to \eqref{eq_classical}}: For any $\varphi\in C(M)$, the function $\Psma\varphi:M\to \R$ is a viscosity solution of the Hamilton-Jacobi equation $G(x, d_x u)=c(G)$.
    	
    	\item [(iii)] \textbf{Fixed Point Characterization}: A continuous function $u:M\to \R$ is a viscosity solution of the Hamilton-Jacobi equation $G(x, d_x u)=c(G)$ if and only if it is a fixed point of the operator $\Psma$, namely,
    	\[
    	u(x)=\inf_{\mu\in \mathfrak{M}(L_G)} \left\{ \frac{\int_{M}  \sigma(y) h_G(y,x)\,d\mu(y)}{\int_{M} \sigma(y) \,d\mu(y)   } +\frac{\int_M  \sigma(y)u(y)\,d\mu(y)}{\int_{M} \sigma(y) \,d\mu(y) } \right\}.
        \]
    \end{itemize}
 \end{THEalpha}
 
 \,
 
Theorem \ref{mainresult3} (ii)-(iii) implies that $\Psma\circ \Psma=\Psma$, and we also directly obtain the following result.
 \begin{CORalpha}\label{mainre31}
 For any fixed continuous positive function $\sigma$, the image of the operator
	$\Psma$ in $C(M)$ coincides with the set of all viscosity solutions to the Hamilton-Jacobi equation $G(x, d_xu)=c(G)$. 
 \end{CORalpha}

\begin{Rem}
Thus, the operator is indeed a solution operator, and it gives a novel variational construction of viscosity solutions (or weak KAM solutions), providing a useful insight into classical weak KAM theory and Aubry-Mather theory.
\end{Rem}

\subsubsection{Uniqueness sets for equation \texorpdfstring{\eqref{eq_classical}}{}}
The Hamilton-Jacobi equation \eqref{eq_classical} admits infinitely many solutions. According to \cite{Fathi_Siconolfi2005}, the projected Mather set $\mathcal{M}_{L_G}$ (see Definition \ref{Def_Matherset}) serves as a {\it uniqueness set}: if two solutions $w_1$ and $w_2$ satisfy $w_1\leq w_2$ on $\mathcal{M}_{L_G}$, then $w_1\leq w_2$ on $M$. The connection between the operator $\Psma$ and viscosity solutions, as shown in  Theorem \ref{mainresult3}, allows to refine the classical uniqueness result as follows:

\begin{THEalpha}\label{mainresult4}
Assume conditions {\bf(G1)}--{\bf(G2)}. Let $\sigma: M\to (0,+\infty)$ be a continuous function.  If $u_1$ and $u_2$ are viscosity solutions to the Hamilton-Jacobi equation \eqref{eq_classical} and 
\begin{equation*}
	\int_M \sigma(y) u_1(y)\,d\mu\leq \int_M \sigma(y) u_2(y)\,d\mu
\end{equation*}
for all projected Mather measures $\mu$ of $L_G$, then $u_1\leq u_2$ on $M$.
\end{THEalpha}

This provides a comparison principle based on integrals with respect to Mather measures.

\begin{CORalpha}\label{mainresult5}
In the special case $\sigma(x)\equiv 1$ on $M$, Theorem \ref{mainresult4} simplifies to the following:
for any viscosity solutions $u_1$ and $u_2$  of equation \eqref{eq_classical}, if   
\begin{equation*}
	\int_M u_1(y)\,d\mu\leq \int_M u_2(y)\,d\mu
\end{equation*}
for all projected Mather measures $\mu$ of $L_G$, then $u_1\leq u_2$ on $M$.
\end{CORalpha}

The uniqueness result in Corollary \ref{mainresult5} aligns with an earlier result of  Mitake and Tran \cite{Mitake_Tran_uniqueness_set}, who used a second-order PDE approach with an additional assumption:
\begin{equation*}
	\lim_{\|p\|_x\to \infty}\frac{1}{2} \big(G(x,p))^2+\partial_x G(x,p) \cdot p=+\infty.
\end{equation*}
In contrast, our approach, using tools from Aubry-Mather and weak KAM theory, does not need such a growth condition. However, we remark that the approach of \cite{Mitake_Tran_uniqueness_set} is also applicable in studying the non-uniqueness phenomenon for certain second order (degenerate viscous) H-J equations.

\subsubsection{Equilibrium Mather measure}

Let us now focus on the operator $\Psma$ for the special case $\sigma\equiv 1$, denoted as $\Psmone$ for brevity.
We are interested in measures that  realize the infimum in the definition of $\Psmone\varphi(x)$.

Consider a function $\varphi\in C(M)$ and a point $x\in M$. In this work, a Mather measure $\tilde\mu\in \tilde{\mathfrak{M}}(L_G)$ on $TM$ is called 
 an {\em equilibrium Mather measure} for $\varphi$ and $x$ if it satisfies:  
\[ \int_{TM}  h_G(y, x)\,d\tmu(y,v) +\int_{TM} \varphi(y)\,d\tmu(y,v)=\Psmone\varphi(x).\]
In this case, the projection $\mu=\pi_{\#}\tmu$ of $\tmu$ on $M$ satisfies
\[ \int_{M}  h_G(y, x)\,d\mu(y) +\int_{M} \varphi(y)\,d\mu(y)=\Psmone\varphi(x).\]
Such a $\mu$ is referred to as the {\em equilibrium projected Mather measure} for $\varphi$ and $x$.

Uniqueness of equilibrium Mather measures is a delicate problem. To explore it, we denote by $\tilde{\mathcal{E}_{G}}(\varphi, x)$ the set of all equilibrium Mather measures for $\varphi$ and $x$, and by $\mathcal{E}_{G}(\varphi, x)$ the set of all equilibrium projected Mather measures for $\varphi$ and $x$. A natural question arises: 

\begin{Prbm}
{\em When does the set $\tilde{\mathcal{E}_{G}}(\varphi, x)$ or $\mathcal{E}_{G}(\varphi, x)$ contain just one  element?}
\end{Prbm}

This of course happens when the Lagrangian $L_G$ has just a unique Mather measure. In the case where the Lagrangian $L_G$ has a lot of Mather measures, i.e., $\operatorname{card}(\tilde{\mathfrak{M}}(L_G))> 1$, we have the following structural property: 

\begin{THEalpha}\label{mainresult6}
	Let $G: T^*M\to \R $ be a $C^2$ Tonelli Hamiltonian. Then
	\begin{enumerate}
		\item [(i)] For any $\varphi\in C(M)$ and $x\in M$, the set $\tilde{\mathcal{E}_{G}}(\varphi, x)$ is convex and compact. Moreover, $\tilde{\mathcal{E}_{G}}(\varphi, x)$ contains ergodic invariant measures: the extremal points of the convex set $\tilde{\mathcal{E}_{G}}(\varphi, x)$ are ergodic. 
	
	\item [(ii)] Given $x\in M$. There is a  residual \footnote{A set is called residual if it contains a countable intersection of open and dense subsets.} subset $\mathcal{O}_x\subset C(M)$ such that for any $\varphi\in \mathcal{O}_x$, 
	\[\operatorname{card}(\tilde{\mathcal{E}_{G}}(\varphi, x))=\operatorname{card}(\mathcal{E}_{G}(\varphi, x))=1.\]
	\item  [(iii)] There is a dense subset $\mathcal{D}\subset M$ and a residual subset $\mathcal{O}\subset C(M)$ such that for any $z\in \mathcal{D}$ and $\varphi\in \mathcal{O}$,
		\[\operatorname{card}(\tilde{\mathcal{E}_{G}}(\varphi, z))=\operatorname{card}(\mathcal{E}_{G}(\varphi, z))=1.\]
	\end{enumerate}
\end{THEalpha}

Recall that the term $C^2$  Tonelli Hamiltonian refers to a Hamiltonian $G$ satisfying the superlinearity condition {\bf (G2)}, and strict convexity in the fibers in the sense that for every $(x,p)\in T^*M$ the second derivative along the fibers $\partial^2H/\partial p^2(x,p)$ is positive definite.

\,

\noindent \textbf{Organization of the paper}. This paper is organized as follows. Section \ref{section_UniformEstimates} is devoted to establishing uniform estimates for the solutions of the perturbed equations, while Section \ref{section_calibratingcurves} discusses the properties of $u_\lambda$-calibrated curves. In Section \ref{section_approximation_Mathermeasures}, we explore the approximation of Mather measures through measures supported on the calibrated curves, and in Section \ref{section_asymp_behavior}, we analyze the asymptotic behavior of the solutions as $\lambda$ tends to zero, thereby proving Theorem \ref{mainresult1}. Section \ref{section_Anotherformula} provides an alternative formula for the limit solution (Theorem \ref{mainresult2}). Finally, Section \ref{section_proofofCandD} is devoted to the proof of Theorem \ref{mainresult3}, Theorem \ref{mainresult4} and Theorem \ref{mainresult6}. The Appendix includes a review of weak KAM theory and Aubry-Mather theory for convex Hamiltonians and Lagrangians, which form the theoretical foundation for our analysis.

%%%%%%%%%%%%%%%%%%%%%%%%%%%%%%%%%%%%%%%%%%%%%%%%%%%%%%%%%%%%%%%%%%

\section{Uniform estimates for the solutions of equation \texorpdfstring{\eqref{eq_H_lambda}}{}}\label{section_UniformEstimates}

Throughout this section, we always assume that conditions {\bf(H1)--\bf(H4)} are satisfied.

We first observe that the critical Hamilton-Jacobi equation $H(x, d_xu,0)=c(H^0)$ admits continuous viscosity solutions, which are uniformly Lipschitz continuous on $M$ with a common Lipschitz constant $K_0=\sup\{\|p\|_x: H(x,p,0)\leq c(H^0) \}.$ Let us define 
\[S_0:=\{(x,p)\in T^*M : x\in M, \|p\|_x\leq K_0\},\] which is a compact subset of $T^*M$. By condition {\bf(H4)}, $\frac{\partial H}{\partial u}(x,p,0)>0$ on $T^*M$, ensuring the existence of a constant $\delta=\delta(S_0)>0$ such that 
\begin{equation}\label{tks0} 
\delta\leq \frac{\partial H}{\partial u}(x,p,0), \quad \forall~(x,p)\in S_0.
\end{equation}

Since $V(x,\lambda)$ converges uniformly to $V(x,0)$ as $\lambda\to 0$ and $M$ is compact, we can choose a constant $\mfK> 1$ such that 
\begin{equation}\label{cka}
	|V(x,\lambda)|\leq \mfK-1
\end{equation}
for all $x\in M$ and sufficiently small $\lambda$.

Now, let $\hat u:M\to \R$ be a viscosity solution of the critical Hamilton-Jacobi equation $H(x, d_xu,0)$ $=$ $c(H^0)$. We then define two functions $\hat{u}_-$ and $\hat{u}_+$ on $M$ as follows:
\begin{equation}\label{defupm}
\hat{u}_-=\hat{u}-\max_{M} \hat{u}-\mfK\delta^{-1}, \qquad \hat{u}_+=\hat{u}-\min_{M} \hat{u}+\mfK\delta^{-1}.	
\end{equation}
Notably, we have
\begin{equation}\label{fnj2}
 \|\hat{u}_-\|_\infty\leq  T,\qquad \|\hat{u}_+\|_\infty\leq T,	
\end{equation} 
where $T= K_0\, \textup{diam}(M)+\mfK\delta^{-1}$. 

Then, the following result holds.

\begin{Lem}\label{lem_subsup0}
	 We can find $\lambda_0>0$, such that the functions $\hat{u}_-:M\to\R $ and $\hat{u}_+:M\to\R$ defined in \eqref{defupm} are, respectively, viscosity subsolution and supersolution of equation \eqref{eq_H_lambda} for all $\lambda\in (0, \lambda_0]$.
\end{Lem}
\begin{proof}
The proof relies on condition {\bf(H4)}, which provides a modulus of continuity $\eta_{_{S_0}}$ such that
\begin{equation}\label{cka10}
\left|H(x,p,u)-H(x,p,0)-u\frac{\partial H}{\partial u}(x,p,0)\right|\leq |u|\eta_{_{S_0}}(|u|),\quad\forall~(x,p,u)\in S_0\times\R.
\end{equation}

We now show that $\hat{u}_-$ is a viscosity subsolution of equation \eqref{eq_H_lambda} with  $\lambda>0$ being suitably small. This involves showing that 
\begin{equation}\label{wtst1}
H(x, p', \lambda \hat{u}_-(x))+\lambda V(x,\lambda) \leq c(H^0),\quad\text{$\forall~ x\in M$ and $p'\in D^+\hat{u}_-(x)$},
\end{equation}
where $D^+\hat{u}_-(x)$ denotes the superdifferential of $\hat{u}_-$  at $x$, see \cite{Barles_book, Bardi_Capuzzo-Dolcetta1997}. Note that $D^+ \hat{u}_-(x)=D^+ \hat{u}(x)$, so for every $p'\in D^+ \hat{u}_-(x)$ we have \[H(x, p',0)\leq c(H^0),\qquad \|p'\|_x\leq K_0.\]
By \eqref{cka10} and $\hat{u}_-\leq -\mfK\delta^{-1}$, 
    \begin{align*}
	 H(x, p', \lambda \hat{u}_-(x)) \leq & H(x, p',0)+\lambda \hat{u}_-(x)\frac{\partial H}{\partial u}(x,p',0)+\| \lambda \hat{u}_-\|_\infty\, \eta_{_{S_0}}(\|\lambda \hat{u}_-\|_\infty)\\
	  \leq & c(H^0)-\lambda \mfK\delta^{-1}\frac{\partial H}{\partial u}(x,p',0) +\lambda\|\hat{u}_-\|_\infty\,\eta_{_{S_0}}(\lambda\|\hat{u}_-\|_\infty). 
	 \end{align*}
Using inequalities \eqref{tks0}, \eqref{cka} and \eqref{fnj2} we get
\begin{align*}
	H(x, p', \lambda \hat{u}_-(x))+\lambda V(x,\lambda) \leq & c(H^0)-\lambda \mfK\delta^{-1}\delta +\lambda T\,\eta_{_{S_0}}(\lambda T)+\lambda(\mfK-1)  \\
	 = & c(H^0)+ \lambda\left(T\,\eta_{_{S_0}}(\lambda T)-1\right).
\end{align*}
We will show below that 
\begin{equation}\label{aiwh3}
T\eta_{_{S_0}}\left(\lambda T\right)\leq 1
\end{equation}
when $\lambda$ is sufficiently small.  Recall that  $\eta_{_{S_0}}:[0,\infty)\to [0,\infty)$ is non-decreasing and $\eta_{_{S_0}}(r)\to 0$ as $r$ tends to zero. Thus we can find a small constant $\lambda_0>0$, such that 
\begin{equation}\label{nxzq}
	\eta_{_{S_0}}\left(\lambda T\right)\leq T^{-1},\quad \text{for all $\lambda\in(0,\lambda_0]$.}
\end{equation}
This verifies \eqref{aiwh3} for all $\lambda\in(0,\lambda_0]$. Therefore, we can conclude that \eqref{wtst1} holds for all $\lambda\in(0,\lambda_0]$, so $\hat{u}_-$ is a viscosity subsolution of equation \eqref{eq_H_lambda}.

To show that $\hat{u}_+$ is a viscosity supersolution of equation \eqref{eq_H_lambda} for  $\lambda\in (0,\lambda_0]$, we only need to check that  
\begin{equation}\label{wtst2}
H(x, p'', \lambda \hat{u}_+(x))+\lambda V(x,\lambda) \geq c(H^0),\quad\text{$\forall~x\in M$ and $p''\in D^-\hat{u}_+(x)$},
\end{equation}
where $D^-\hat{u}_+(x)$ denotes the subdifferential of $\hat{u}_+$ at $x$. Note that $D^- \hat{u}_+(x)=D^- \hat{u}(x)$ and the Lipschitz constant $\textup{Lip}(\hat{u}_+)\leq K_0$, then for every $p''\in D^- \hat{u}_+(x)$ we have
 \[\|p''\|_x\leq K_0,\qquad H(x, p'',0)\geq c(H^0).\] 
By \eqref{cka10}, together with $\hat{u}_+\geq \mfK\delta^{-1}$ and $\|\hat{u}_+\|_\infty\leq T$,  we obtain
\begin{align*}
	 H(x, p'', \lambda \hat{u}_+(x)) \geq & H(x, p'',0)+\lambda \hat{u}_+(x)\frac{\partial H}{\partial u}(x,p'',0)-\| \lambda \hat{u}_+\|_\infty\, \eta_{_{S_0}}(\|\lambda \hat{u}_+\|_\infty)\\
              \geq & c(H^0)+\lambda \mfK\delta^{-1}\frac{\partial H}{\partial u}(x,p'',0) -\lambda T \eta_{_{S_0}}(\lambda T)\\ 
              \geq  & c(H^0)+\lambda \mfK -\lambda 
\end{align*}
where the last line follows from \eqref{tks0} and \eqref{nxzq}. In view of \eqref{cka}, we conclude that
\begin{align*}
	H(x, p'', \lambda \hat{u}_+(x))+\lambda V(x,\lambda) 
	 \geq & c(H^0)+\lambda \mfK -\lambda -\lambda (\mfK-1)=c(H_0)
\end{align*}
for all $\lambda\in (0,\lambda_0]$. This therefore verifies \eqref{wtst2}, so $\hat{u}_+$ is a viscosity supersolution of equation \eqref{eq_H_lambda} for $0<\lambda\leq \lambda_0$.

\end{proof}

Based on Lemma \ref{lem_subsup0} and $\hat u_-\leq \hat u_+$, we can utilize the Perron method \cite{Ishii1987} to construct a possible solution to \eqref{eq_H_lambda}, $\lambda\in (0,\lambda_0]$. In fact, the function defined by
\begin{equation}\label{constulam} 
u_\lambda(x):=\sup\{w(x) \mid \hat u_-\leq w\leq \hat u_+, w \text{~is a continuous viscosity subsolution of~} \eqref{eq_H_lambda} \}
\end{equation}	 
is indeed a continuous viscosity solution on $M$. Besides, by the monotonicity condition {\bf(H3)}, equation \eqref{eq_H_lambda} obeys the comparison principle, so $u_\lambda:M\to\R$ is the unique continuous viscosity solution. In conclusion, the following proposition is derived:

\begin{Pro}\label{Thm_existencesolutions}
There exists a constant $\lambda_0>0$ such that for each $\lambda\in (0, \lambda_0]$, equation \eqref{eq_H_lambda} has a unique continuous viscosity solution $u_\lam: M\to \R$. 

Furthermore, the family $(u_\lambda)_{\lambda\in (0,\lambda_0]}$ is equibounded and equi-Lipschitz.  
\end{Pro}

\begin{proof}
The first part has been proved above. By the construction \eqref{constulam}, it is clear that 
\begin{equation*}
	\|u_\lambda\|_\infty \leq  \max\{ \|\hat{u}_-\|_\infty, \|\hat{u}_+\|_\infty\} \leq T=K_0\, \textup{diam}(M)+\mfK\delta^{-1},
\end{equation*}
so the family $\{u_\lambda: 0<\lambda\leq\lambda_0\}$ is equibounded. Now, it remains to show the equi-Lipschitz property. Indeed, by the monotonicity condition {\bf(H3)} we easily see that
\[  H(x,p,-\lambda_0 T)\leq H(x,p,\lambda u_\lambda(x)),\quad \text{for all $(x,p)\in T^*M$ and $\lambda\in (0,\lambda_0]$}. \]
Hence, all viscosity solutions $u_\lam : M\to\R$ of \eqref{eq_H_lambda}, $\lambda\in (0,\lambda_0]$, are viscosity subsolutions of the  equation
\[  H(x, d_x u_\lambda, -\lambda_0 T)+\lambda V(x,\lambda)= c(H^0). \] 
Together with \eqref{cka}, $(u_\lam)_{\lambda\in (0,\lambda_0]}$ are thus viscosity subsolutions of the following equation
\[H^{-\lambda_0 T}(x, d_xu_\lambda)=c(H^0)+\lambda_0 (\mfK-1), \]
where the Hamiltonian $H^{-\lambda_0 T}(x, p)=H(x,p,-\lambda_0 T)$. 
Since $H^{-\lambda_0 T}$ satisfies the superlinearity condition, the family $(u_\lam)_{\lambda\in (0,\lambda_0]}$ is equi-Lipschitz with  Lipschitz constant 
$\max\{\|p\|_x : H^{-\lambda_0 T}(x,p)\leq c(H^0)+\lambda_0 (\mfK-1)\}$. This finishes the proof.
\end{proof}

We stress that Proposition \ref{Thm_existencesolutions} may not hold for $\lambda>\lambda_0$. We provide below a concrete example to illustrate this phenomenon. 

\begin{Ex}[Non-existence of solutions for large $\lambda$]\label{ex_hasnosol}
Consider a Hamiltonian 
\[H(x,p,u)=f(u)+\|p\|_x^2,  \quad \forall ~ (x,p,u)\in T^*M\times\R\] with $f(u)=\arctan u$. It is straightforward to verify that $H$ satisfies conditions {\bf(H1)-(H4)}, and the critical value $c(H^0)=0$. Let $V(x,\lambda)\equiv \frac{\pi}{2}$, and consider the equation:
	\begin{equation}\label{qeqdniwq1}
f(\lambda u_\lambda(x))+\|d_x u_\lambda\|_x^2+\lambda V(x,\lambda)=0\quad \text{in $M$.}
	\end{equation}
We claim that for all $\lambda>1$, equation \eqref{qeqdniwq1} has no solutions. Suppose there exists $\lambda_1>1$ such that this equation admits a viscosity solution $u_{\lambda_1}:M\to \R$. Since $M$ is compact, there must exist a maximum point $x_0\in M$ where $u_{\lambda_1}(x_0)=\max_{x\in M} u_{\lambda_1}(x)$. At this point, $0\in D^+u_{\lambda_1}(x_0)$, and by the definition of viscosity subsolutions, we have
\begin{equation}\label{qeqdniwq2}
f(\lambda_1 u_{\lambda_1}(x_0))+\|0\|_{x_0}^2+\lambda_1 V(x_0,\lambda_1)\leq 0.
\end{equation}
Given that $f\geq -\frac{\pi}{2}$ and $V\equiv \frac{\pi}{2}$, this implies  $f(\lambda_1 u_{\lambda_1}(x_0))+\lambda_1 V(x_0,\lambda_1)>0$, which contradicts \eqref{qeqdniwq2}. Therefore, no solutions exist for $\lambda>1$.
\end{Ex}

To explore the dynamical behavior of the viscosity solution $u_\lambda$ of the perturbation problem \eqref{eq_H_lambda}, we can interpret $u_\lambda$ as a solution of the Hamilton-Jacobi equation of a $u$-independent Hamiltonian, thereby allowing us to apply the results from Appendix \ref{sect_weakKAMandAubrymather}. To do this it is convenient to introduce the following notation: consider a continuous function $\psi: M\to \R$, we denote by $\HH^\psi:T^*M\to\R$ a new Hamiltonian given by
\begin{equation}\label{HHdef}
\HH^\psi(x,p)=H(x,p, \psi(x)),\quad \text{for all $(x,p)\in T^*M$}.	
\end{equation}
Since $H$ satisfies the convexity condition {\bf(H1)} and the superlinearity condition {\bf(H2)},  $\HH^\psi$ is also convex and superlinear in the momentum variable $p$.
Accordingly, $\HH^\psi$ has a conjugated Lagrangian which we denote by $\LL^\psi: TM\to\R$,
\begin{equation}\label{LLdef}
	\LL^\psi(x,v)=L(x,v,\psi(x)),\quad \text{for all $(x,v)\in TM$}.	
\end{equation}

Introducing this notation allows us to establish the following: 

\begin{Pro}\label{Pro_newsysproperty}
For $\lambda\in(0,\lambda_0]$, the solution $u_\lambda:M\to\R$ of the perturbation problem \eqref{eq_H_lambda} is also a viscosity solution of the equation: 	
\begin{equation*}
		\HH^{\lambda u_\lambda}(x, d_x u)+\lambda V(x,\lambda)=c(H^0).
	\end{equation*}
\end{Pro}
\begin{proof}
	This follows directly from the expression that $\HH^{\lambda u_\lambda}(x,p)=H(x,p, \lambda u_\lambda(x))$.
\end{proof}

Proposition \ref{Pro_newsysproperty} also implies the following result.
\begin{Cor}\label{Cor_critvalue}
The critical value of the Hamiltonian $\HH^{\lambda u_\lambda}+\lambda V$ is exactly $c(H^0)$. 
\end{Cor}

Since $\HH^{\lambda u_\lambda}+\lambda V$ has a conjugated Lagrangian $\LL^{\lambda u_\lambda}-\lambda V$, the following result is derived.
\begin{Pro}\label{Pro_clomeasgeq0}
	For the solution $u_\lambda:M\to\R$ of the perturbation problem \eqref{eq_H_lambda} with $0<\lambda\leq\lambda_0$, we have 	
	\begin{equation*}
		\int_{TM} \left[L(x,v,\lambda u_\lambda(x))-\lambda V(x,\lambda)+c(H^0)\right]\,d\tilde\mu\geq 0, \quad \text{for any closed measure $\tilde\mu$ on $TM$.}
    \end{equation*}
\end{Pro}
\begin{proof}
 Due to Corollary \ref{Cor_critvalue}, we infer from Proposition \ref{Propertiesclosedmeas} that for any closed measure $\tilde\mu$ on $TM$,
	\begin{equation*}
		\int_{TM} \left[\LL^{\lambda u_\lambda}(x,v)-\lambda V(x,\lambda)\right]\,d\tilde\mu\geq -c(H^0).
    \end{equation*}
    Given that $\LL^{\lambda u_\lambda}(x,v)=L(x,v,\lambda u_\lambda(x))$ by \eqref{LLdef}, it follows that
    \[\int_{TM} \left[L(x,v,\lambda u_\lambda(x))-\lambda V(x,\lambda)+c(H^0)\right]\,d\tilde\mu\geq 0.\]
    This finishes the proof.
    \end{proof}

\section{\texorpdfstring{$u_\lambda$}{}-calibrated curves}\label{section_calibratingcurves}

This section discusses the properties of $u_\lambda$-calibrated curves, which are important to understand the structure and asymptotic behavior of the solution $u_\lambda$ of  \eqref{eq_H_lambda}. These curves provide insights into how solutions evolve and help establish connections between the Hamiltonian dynamics and the corresponding Lagrangian formulation.

Owing to Proposition \ref{Pro_newsysproperty}, we derive the following result.
\begin{Pro}\label{Pro_solu_cali}
For the viscosity solution $u_\lambda:M\to \R$ of \eqref{eq_H_lambda}, with $\lambda\in (0,\lambda_0]$, the following statements hold:
	\begin{enumerate}
		\item[(1)] For every absolutely continuous curve $\xi:[a,b]\to M$, 
\begin{equation}\label{dominatedcurves}
u_\lambda(\xi(b))-u_\lambda(\xi(a))\leq 
\int_a^b \left [L\big(\xi(s),\dot\xi(s), \lambda u_\lambda(\xi(s))\big)-\lambda V(\xi(s),\lambda) +c(H^0)\right] \,ds.
\end{equation}
\item[(2)] For each $x\in M$, there exists an absolutely continuous curve $\xilax:(-\infty,0]\to M$ with $\xilax(0)=x$, such that for any $t\geq 0$,
\begin{equation}\label{calibratedcurves} 
u_\lambda(\xilax(0))-u_\lambda(\xilax(-t))=
\int_{-t}^0 \left [ L\big(\xilax(s),\dot\xilax(s), \lambda  u_\lambda(\xilax(s))\big)-\lambda V(\xilax(s),\lambda)+c(H^0) \right ] \,ds.
\end{equation}
	\end{enumerate}
\end{Pro}
\begin{proof}
Since $u_\lambda$ is a viscosity solution of the equation 
$\HH^{\lambda u_\lambda}(x, d_x u)+\lambda V(x,\lambda)=c(H^0),$  
and the Lagrangian associated to $\HH^{\lambda u_\lambda}+\lambda V$ is precisely $\LL^{\lambda u_\lambda}-\lambda V$,
it is thus classical that for every absolutely continuous curve $\xi:[a,b]\to M$, 
\begin{equation*}
u_\lambda(\xi(b))-u_\lambda(\xi(a))\leq 
\int_a^b\left[ \LL^{\lambda u_\lambda}(\xi(s),\dot\xi(s))-\lambda V(\xi(s),\lambda) +c(H^0)\right] \,ds.
\end{equation*}
In addition, we also have that for each $x\in M$, there exists an absolutely continuous curve $\xilax:(-\infty,0]\to M$ with $\xilax(0)=x$, such that for any $t\geq 0$,
\begin{equation*}
u_\lambda(\xilax(0))-u_\lambda(\xilax(-t))=
\int_{-t}^0 \left[\LL^{\lambda u_\lambda}(\xilax(s),\dot\xilax(s))-\lambda V(\xilax(s),\lambda)+c(H^0)\right]\,ds.
\end{equation*}
Given that $\LL^{\lambda u_\lambda}(x,v)=L(x,v,\lambda u_\lambda(x))$,  the statements \eqref{dominatedcurves} and \eqref{calibratedcurves} are verified.
\end{proof}

\begin{defn}
A curve $\gamma:I\to M$, where $I$ is an interval, is called {\em $u_\lambda$-calibrated} if for any subinterval $[a,b]\subset I$, 
\begin{equation*}
u_\lambda(\gamma(b))-u_\lambda(\gamma(a))=
\int_a^b \left[ L\big(\gamma(s),\dot\gamma(s), \lambda u_\lambda(\gamma(s))\big)-\lambda V(\gamma(s),\lambda) +c(H^0)\right] \,ds.
\end{equation*}
\end{defn}

By the superlinearity of $L$, it is known that calibrated curves are Lipschitz continuous. The following result shows that the $u_\lambda$-calibrated curves are indeed equi-Lipschitz for $0<\lambda\leq\lambda_0$. 
\begin{Pro}\label{Pro_calibratingLip}
There exists a constant $\mathfrak{B}=\mathfrak{B}(\lambda_0)>0$, such that for any curve $\gamma:[a,b]\to M$ which is $u_\lambda$-calibrated, $\lambda\in (0,\lambda_0]$, $\gamma$ is Lipschitz continuous with Lipschitz constant $\mathfrak{B}$.
\end{Pro}

\begin{proof}
According to Proposition \ref{Thm_existencesolutions}, we can find a constant $\mfD_0=\mfD_0(\lambda_0)>0$ such that
\begin{equation}\label{ncoqd}
\|u_\lambda\|_\infty\leq \mfD_0 \quad \textup{and} \quad   \textup{Lip}(u_\lambda)\leq \mfD_0, \quad \text{for all $\lambda\in (0,\lambda_0]$.}
\end{equation}
Consider $\lambda\in (0,\lambda_0]$ and a $u_\lambda$-calibrated curve  $\gamma:[a,b]\to M$. For any subinterval $[c,d]\subset [a,b]$, as $L(x,v,u)$ is strictly decreasing in $u$, we have
\begin{equation}\label{eq_j28en}
\begin{aligned}
	u_\lambda(\gamma(d))-u_\lambda(\gamma(c))=&\int_c^d L\big(\gamma(t),\dot\gamma(t), \lambda  u_\lambda(\gamma(t))\big)-\lambda V(\gamma(t),\lambda)+c(H^0)\,dt\\ 
	\geq & \int_c^d L(\gamma(t),\dot\gamma(t), \lambda_0  \mfD_0)-\lambda V(\gamma(t),\lambda)+c(H^0)\,dt.
\end{aligned}
\end{equation}
Recall that by \eqref{cka}, $| V(\cdot,\lambda)|\leq (\mfK-1)$, and by the superlinearity condition of $L$, there exists a constant $C_0>0$ such that 
\begin{equation}\label{supajdn}
L(x,v,\lambda_0 \mfD_0)\geq (\mfD_0+1) \|v\|_x-C_0, \quad\text{for all $(x,v)\in TM$}.
\end{equation}
Thus, we deduce from \eqref{eq_j28en} and \eqref{supajdn} that 
\begin{equation}\label{eq_loew62}
\begin{aligned}
		u_\lambda(\gamma(d))-u_\lambda(\gamma(c))\geq &\int_c^d (\mfD_0+1) \|\dot\gamma(t)\|_{\gamma(t)} -C_0-\lambda\,(\mfK-1) +c(H^0)\, dt\\
\geq &(\mfD_0+1) \textup{dist}(\gamma(c),\gamma(d))-(d-c)\big[C_0+\lambda_0\,(\mfK-1)-c(H^0)\big].
\end{aligned}
\end{equation}
By \eqref{ncoqd} we also have $u_\lambda(\gamma(d))-u_\lambda(\gamma(c))\leq \mfD_0\, \textup{dist}(\gamma(c),\gamma(d))$. Hence, \eqref{eq_loew62} leads to
\[(d-c)\big[C_0+\lambda_0\,(\mfK-1)-c(H^0)\big]\geq \textup{dist}(\gamma(c),\gamma(d)),\]
which in turn implies 
\begin{equation*}
	\frac{\textup{dist}(\gamma(c), \gamma(d))}{d-c}\leq C_0+\lambda_0\,(\mfK-1)-c(H^0).
\end{equation*}
As this inequality holds for any $[c,d]\subset [a,b]$, we conclude that $\gamma:[a,b]\to M$
 is Lipschitz with Lipschitz constant $\mathfrak{B}=C_0+\lambda_0\,(\mfK-1)-c(H^0)$.
\end{proof}

As a consequence, the solution $u_\lambda$ exhibits the following differential property along the $u_\lambda$-calibrated curves.
\begin{Cor}\label{Cor_differulam}
	Let $\gamma:[a,b]\to M$ be any $u_\lambda$-calibrated curve, $0<\lambda\leq \lambda_0$, then the function  $t\mapsto u_\lambda( \gamma(t))$ is differentiable almost everywhere on $[a,b]$. Moreover, for every $t\in(a,b)$ where $u_\lambda(\gamma(\cdot))$ is differentiable, we have
\begin{equation}\label{dulambdadt}
	\frac{d}{dt} u_\lambda(\gamma(t))=L\big(\gamma(t),\dot\gamma(t), \lambda  u_\lambda(\gamma(t))\big)-\lambda V(\gamma(t),\lambda)+c(H^0).
\end{equation}
\end{Cor}
\begin{proof}
The function $u_\lambda(\gamma(t))$ is the composition of the Lipschitz function $u_\lambda$ with the mapping $t\mapsto \gamma(t)$. By Proposition \ref{Pro_calibratingLip}, the  $u_\lambda$-calibrated curve $\gamma(t)$ is  Lipschitz, so the function  $t\mapsto u_\lambda( \gamma(t))$ is also Lipschitz continuous. Thus, it is differentiable almost everywhere on $[a,b]$. 

Since $\gamma|_{[a,b]}$ is $u_\lambda$-calibrated, it is well known that for every $t\in (a,b)$, the restriction $\gamma|_{[a,t]}$ is also $u_\lambda$-calibrated. Hence,
\begin{equation*}
	u_\lambda(\gamma(t))-u_\lambda(\gamma(a))=
\int_a^t \left[L\big(\gamma(s),\dot\gamma(s), \lambda u_\lambda(\gamma(s))\big)-\lambda V(\gamma(s),\lambda) +c(H^0)\right] \,ds.
\end{equation*}
Therefore, for each $t$ where $u_\lambda( \gamma(\cdot))$ is differentiable, differentiating the above equation with respect to $t$ yields the desired equality \eqref{dulambdadt}.
\end{proof}

%%%%%%%%%%%%%%%%%%%%%%%%%%%%%%%%%%%%%%%%%%%%%%%%%%%%%%%%%%%%%%%%%%
\section{Approximation of Mather Measures}\label{section_approximation_Mathermeasures}
Mather measures, named after J. Mather, are probability measures that minimize the action of the Lagrangian system. See Appendix \ref{sect_weakKAMandAubrymather} for more detail.  They are critical in understanding the long-term behavior of dynamical systems. This section is concerned with the approximation of Mather measures, particularly focusing on the approximation by measures having their supports on the $u_\lambda$-calibrated curves. The results in this section demonstrate how Mather measures can be approximated using solutions to perturbed Hamilton-Jacobi equations. Such an approximation is crucial for understanding the asymptotic behavior of the perturbed solution $u_\lambda$ as $\lambda$ tends to zero.

\subsection{Probability measures supported on calibrated curves}

Let $u_\lambda:M\to\R$ be the viscosity solution of equation \eqref{eq_H_lambda} with $0<\lambda\leq\lambda_0$, where $\lambda_0$ is given by Proposition \ref{Thm_existencesolutions}. For any $x\in M$, we can find a backward $u_\lambda$-calibrated curve $\xilax:(-\infty, 0]\to M$ with $\xilax(0)=x$, as established in Proposition \ref{Pro_solu_cali}. 
\begin{defn}\label{def_measoncalicurve}
Let $x\in M$ and $\lambda\in (0,\lambda_0]$. For the $u_\lambda$-calibrated curve $\xilax:(-\infty, 0]\to M$, we define a probability measure $\tmulax=\tilde\mu^\lambda_{\xilax}$ on $TM$ by setting
	\begin{equation}\label{def_pro_meas}
	\begin{aligned}
	\int_{TM}  f(y, & v)  \,d\tmulax(y,v):=\\
	 &\frac{1}{ \int_{-\infty}^0   e^{\lam\int^0_t \frac{\partial L}{\partial u }(\xilax(s),\dxilax(s),0)\, ds}\, dt}\int_{-\infty}^0    f(\xilax(t),\dxilax(t))\,  e^{\lam\int^0_t \frac{\partial L}{\partial u }(\xilax(s),\dxilax(s),0)\, ds} \, dt,
	\end{aligned}
	\end{equation}
for any function $f\in C_c(TM)$.
\end{defn}

\begin{Rem}
The  measure $\tmulax$ is well defined: indeed, 	since $\frac{\partial L}{\partial u}(x,v,0)$ is continuous on $TM$ and the curve $\xilax$ is $\mathfrak{B}$-Lipschitz  (as per Proposition \ref{Pro_calibratingLip}), we can find two constants $\delta_1 \geq \delta_2>0$ such that for any $t\in (-\infty,0]$, 
\begin{equation*}
	\delta_1 t\leq \int^0_t \frac{\partial L}{\partial u}(\xilax(s),\dxilax(s),0)\,ds\leq \delta_2 t,
\end{equation*}
which yields 
\begin{equation}\label{na5sj17}
	e^{\lambda\delta_1 t}\leq e^{\lambda\int^0_t \frac{\partial L}{\partial u}(\xilax(s),\dxilax(s),0)\,ds}\leq e^{\lambda\delta_2 t}.
\end{equation}
This therefore ensures that the integrals involved in the definition of $\tmulax$ are finite..  
\end{Rem}

As a result of \eqref{na5sj17}, we directly obtain the following: 
\begin{Lem}\label{dw72je}
	For every backward $u_\lambda$-calibrated curve $\xilax:(-\infty, 0]\to M$, we have
	\begin{equation*}
	 e^{\lambda\int^0_{-\infty} \frac{\partial L}{\partial u}(\xilax(s),\dxilax(s),0)\,ds}=0.
\end{equation*}
\end{Lem}

Now, consider $f(x,v)=\frac{\partial L}{\partial u}(x,v,0)$ in \eqref{def_pro_meas}, we have a lemma that will be useful later.

\begin{Lem}\label{Lem_measidentity}
For the probability measure $\tmulax$ associated with the $u_\lambda$-calibrated curve $\xilax$, 
	\begin{equation*}
		\int_{TM} \frac{\partial L}{\partial  u }(y,v,0)\,d\tmulax(y,v)=-\lambda^{-1}\left( \int_{-\infty}^0   e^{\lam\int^0_t\frac{\partial L}{\partial  u }(\xilax(s),\dxilax(s),0)\, ds}\, dt\right)^{-1}.
	\end{equation*}
\end{Lem}
\begin{proof}
We start by considering the function $f(x,v)=\frac{\partial L}{\partial u}(x,v,0)$ in the definition (see \eqref{def_pro_meas}) of the measure $\tmulax$, which gives
\begin{equation}\label{daoww2}
	\int_{TM} \frac{\partial L}{\partial u }(y,v,0)\,d\tmulax(y,v)
 	=\frac{\int_{-\infty}^0 \frac{\partial L}{\partial  u }(\xilax(t),\dxilax(t),0) \, e^{\lam\int^0_t \frac{\partial L}{ \partial u }(\xilax(s),\dxilax(s),0)\, ds}\, dt}{ \int_{-\infty}^0   e^{\lam\int^0_t \frac{\partial L}{\partial  u }(\xilax(s),\dxilax(s),0)\, ds}\, dt}. 
\end{equation}
Next, we perform an integration by parts on the integral in the numerator. Define
\[F(t):=e^{\lam\int^0_t \frac{\partial L}{\partial u }(\xilax(s),\dxilax(s),0)\, ds},\qquad \forall t\in(-\infty,0].\]
Then, differentiating $F(t)$ with respect to $t$ gives
\begin{equation}\label{useident0} 
	 \frac{d}{dt}F(t)
	 =-\lambda\frac{\partial L}{\partial  u }(\xilax(t),\dxilax(t),0) \, F(t).
\end{equation}
Thus, we can express the integral in the numerator of \eqref{daoww2}  as
\begin{align}\label{hskfnjkqw}
	\int_{-\infty}^0 \frac{\partial L}{\partial  u }(\xilax(t),\dxilax(t),0) \, F(t)\, dt=&-\frac{1}{\lambda}\int_{-\infty}^0   \frac{d}{dt}F(t)\, dt.
\end{align}
Note that by Lemma \ref{dw72je},
\begin{align*}
	\int_{-\infty}^0   \frac{d}{dt}F(t)\, dt=F(0)-\lim_{t\to-\infty} F(t)=1-0=1,
\end{align*}
so the integral \eqref{hskfnjkqw} simplifies to
\begin{equation*}
	\int_{-\infty}^0 \frac{\partial L}{\partial  u }(\xilax(t),\dxilax(t),0) \, F(t)\, dt=-\frac{1}{\lambda}.
\end{equation*}
Substituting this into the expression \eqref{daoww2}, we obtain
\begin{equation*}
	\int_{TM} \frac{\partial L}{\partial u }(y,v,0)\,d\tmulax(y,v)
 	=\frac{-1/\lambda}{ \int_{-\infty}^0   e^{\lam\int^0_t \frac{\partial L}{\partial  u }(\xilax(s),\dxilax(s),0)\, ds}\, dt}. 
\end{equation*}
This completes the proof.
\end{proof}

\subsection{Curve based approximation of Mather measures} By the results from the previous subsection, we now give an approximation scheme for constructing Mather measures. It demonstrates how probability measures supported on $u_\lambda$-calibrated curves can be used to approximate Mather measures of the Lagrangian $L^0(x,v)=L(x,v,0)$.

\begin{Pro}\label{Pro_weakconvergence}
Let $x\in M$. The family of probability measures $\tmulax$, $\lam\in(0,\lam_0]$, is  relatively compact for the weak topology in the space of probability measures on $TM$.

 Moreover, if $\tilde\mu$ is the weak limit of a sequence $\{\tilde{\mu}^{\lambda_k}_x\}$  as $\lambda_k\to 0$, then $\tilde{\mu}$ is  a Mather measure of  the Lagrangian $L^0$, where $L^0(x,v)=L(x,v,0)$. 
\end{Pro}
\begin{proof}
We begin by noting that the measures $\tmulax$, $\lam\in(0,\lam_0]$, are supported on $u_\lambda$-calibrated curves. From Proposition \ref{Pro_calibratingLip}, these curves are equi-Lipschitz with a common Lipschitz constant $\mathfrak{B}$. Thus, the measures $\tmulax$, $\lam\in(0,\lam_0]$ have support contained in a common compact subset of $TM$. As a consequence, the relative compactness of $(\tmulax)_{\lam\in(0,\lam_0]}$ follows. 

Suppose that $\tilde{\mu}^{\lambda_k}_x \rightharpoonup\tilde\mu$ as $k\to \infty$. To show that $\tilde\mu$ is a Mather measure of the Lagrangian $L^0$, we need to verify two key properties: a) $\tilde\mu$ is a closed measure; b) $\tilde\mu$ is minimizing.

\textbf{a) $\tilde\mu$ is a closed measure:}  note that for $\lambda\in (0,\lambda_0]$ and any $C^1$ function $\phi: M\to \R $, by \eqref{def_pro_meas} we directly derive that
\begin{equation}\label{close_mea1}
\begin{aligned}
\int_{TM} d_y\phi(v) \,d\tmulax(y,v)=& \frac{\int_{-\infty}^0 e^{\lam\int^0_t \frac{\partial L}{\partial  u }(\xilax(s),\dxilax(s),0)\,ds}\cdot  \frac{d}{d t}\phi(\xilax(t))\, dt}{\int_{-\infty}^0   e^{\lam\int^0_t \frac{\partial L}{\partial u }(\xilax(s),\dxilax(s),0)\,ds}\,dt}\\
=&\frac{ \phi(x) -\int_{-\infty}^0  \phi(\xilax(t))\, \frac{d}{d t}\left[e^{\lambda\int_{t}^0\frac{\partial L}{\partial u}(\xilax(s),\dxilax(s),0)\,ds} \right]\, dt }{\int_{-\infty}^0   e^{\lam\int^0_t \frac{\partial L}{\partial u }(\xilax(s),\dxilax(s),0)\,ds}\,dt}\\
=&\frac{ \phi(x) +\lambda \int_{-\infty}^0  \phi(\xilax(t))\,\frac{\partial L}{\partial  u }(\xilax(t),\dxilax(t),0) \, e^{\lam\int^0_t \frac{\partial L}{ \partial u }(\xilax(s),\dxilax(s),0)\, ds}\, dt }{\int_{-\infty}^0   e^{\lam\int^0_t \frac{\partial L}{\partial u }(\xilax(s),\dxilax(s),0)\,ds}\,dt}
\end{aligned}
\end{equation}
where the last line has used \eqref{useident0}. Applying Lemma \ref{Lem_measidentity}  and  the definition of $\tmulax$  to \eqref{close_mea1}, we deduce that
\begin{equation}\label{close_shj2} 
	\int_{TM} d_y\phi(v) d\tilde{\mu}^{\lam}_x(y,v)=-\lam \phi(x)\int_{TM} \frac{\partial L}{\partial  u }(y,v,0)d\tmulax(y,v) +\lam \int_{TM}\phi(y) \frac{\partial L}{ \partial u }(y,v,0)d\tmulax(y,v).
\end{equation}
Sending $\lambda\to 0$ in \eqref{close_shj2} yields 
\begin{equation*}
	\lim_{\lam\to 0}\int_{TM} d_y\phi(v) \,d\tilde{\mu}^{\lam}_x(y,v)=0.
\end{equation*}
In particular, if $\tilde{\mu}^{\lambda_k}_x \rightharpoonup\tilde\mu$ as $k\to \infty$, then 
\[ \int_{TM} d_y\phi(v) \,d\tilde{\mu}(y,v)=0.\]
This verifies that $\tilde\mu$ is a closed probability measure.

\textbf{b) $\tilde\mu$ is minimizing:} it suffices to show that 
\begin{equation*}
	\int_{TM} L^0(y,v)+c(H^0)\, d\tilde\mu(y,v)=0.
\end{equation*}
We first observe that for $\lambda\in (0,\lambda_0]$,
\begin{equation}\label{dnw15}
\int_{TM} L^0(y,v)+c(H^0)\,d\tmulax(y,v) =I_1(\lam)+I_2(\lam),	
\end{equation}
where the two summands $I_1$ and $I_2$ are
\begin{align*}
I_1(\lam)=& \int_{TM} L(y,v,\lam u_\lam(y))-\lam V(y,\lam)+c(H^0)\,d\tilde{\mu}_x^\lam(y,v),  \\
I_2(\lam)=&-\int_{TM} L(y,v,\lam u_\lam(y))-\lam V(y,\lam)-L^0(y,v)\,d\tilde{\mu}_x^\lam(y,v).  
\end{align*}
According to Corollary \ref{Cor_differulam}, for each $u_\lam$-calibrated curve $\xilax:(-\infty,0]\to M$, we have
\begin{equation*}
	\frac{d}{dt}u_\lam(\xilax(t))=L\big(\xilax(t),\dxilax(t),\lam u_\lam(\xilax(t))\big)-\lam V(\xilax(t),\lam)+c(H^0),\qquad \textup{a.e. $t\leq0$},
\end{equation*}
which, together with the definition of $\tmulax$, implies that
\begin{align*}
	I_1(\lam)
		&=\frac{\int_{-\infty}^0 \frac{d}{dt} u_\lam(\xilax(t))\cdot e^{\lam\int^0_t \frac{\partial L}{\partial  u }(\xilax(s),\dxilax(s),0)\,ds  }\, dt}{\int_{-\infty}^0    e^{\lam\int^0_t \frac{\partial L}{\partial  u }(\xilax(s),\dxilax(s),0)\,ds  }\, dt}\\
	&=\frac{ u_\lam(x) +\lambda \int_{-\infty}^0 u_\lam(\xilax(t)) \, \frac{\partial L}{\partial  u }(\xilax(t),\dxilax(t),0) \, e^{\lam\int^0_t \frac{\partial L}{\partial  u }(\xilax(s),\dxilax(s),0)\,ds}\, dt }{\int_{-\infty}^0   e^{\lam\int^0_t \frac{\partial L}{\partial u }(\xilax(s),\dxilax(s),0)\,ds}\,dt},
\end{align*}
where the last line has used an integration by parts and \eqref{useident0}.
Then, analogous to the proof of \eqref{close_shj2}, we can use again Lemma \ref{Lem_measidentity} to obtain 
\begin{equation}\label{pjgdi1}
	I_1(\lam)=-\lam u_\lam(x) \int_{TM}  \frac{\partial L}{\partial u }(y,v,0)\, d\tmulax(y,v) +\lam \int_{TM}
u_\lam(y) \frac{\partial L}{ \partial u }(y,v,0)\,d\tmulax(y,v)
\end{equation} 

On the other hand, since the probability measures $\{\tmulax: \lam\in(0,\lam_0]\}$ have support contained in a common compact subset $S=\{(y,v)\in TM : \|v\|_y\leq \mathfrak{B}\}$, we derive from property \eqref{d23rL4} that 
\begin{equation}\label{pjgdi2}
	\big|I_2(\lam)\big|\leq  \left|\int_{TM}\lam u_\lam(y) \frac{\partial L}{\partial u}(y,v,0)+\lam\|u_\lam\|_\infty\eta_{_S}(\lam\|u_\lam\|_\infty)-\lam V(y,\lam) \,d\tilde{\mu}_x^\lam(y,v) \right|.
\end{equation} 
Next, it follows immediately from \eqref{pjgdi1} and \eqref{pjgdi2} that $I_1(\lambda)\to 0$ and $I_2(\lambda)\to 0$ as $\lambda$ tends to zero, so \eqref{dnw15} yields 
\[\lim_{\lam\to 0}\int_{TM} L^0(y,v)+c(H^0)\,d\tmulax(y,v)=0. \]
In particular, if $\tilde{\mu}^{\lambda_k}_x \rightharpoonup\tilde\mu$ as $k\to \infty$, then
\[\int_{TM} L^0(y,v)+c(H^0)\, d\tilde\mu(y,v)=0.\] 
This finally completes the proof.  
\end{proof}

Note that the above result shows that Mather measures can be approximated effectively using solutions of the perturbed Hamilton-Jacobi equations \eqref{eq_H_lambda}.

%%%%%%%%%%%%%%%%%%%%%%%%%%%%%%%%%%%%%%%%%%%%%%%%%%%%%%%%%%%%%%%%%%
\section{Asymptotic behavior of the solutions}\label{section_asymp_behavior}

This section aims to prove Theorem \ref{mainresult1}. As established in Proposition \ref{Thm_existencesolutions}, the family of viscosity solutions $u_\lambda:M\to \R$ of \eqref{eq_H_lambda}, with $\lambda\in (0,\lambda_0]$, is equibounded and equi-Lipschitz continuous.  By the Ascoli-Arzel\`a theorem, the family $(u_\lambda)_{\lambda\in (0,\lambda_0]}$ is relatively compact. Consequently, due to the stability property of viscosity solutions, any convergent subsequence converges to a viscosity solution of the critical equation \eqref{eq_G}. In this section, we will show that all different subsequences yield the same limit.

\subsection{Asymptotic analysis}

In the sequel, $\wtM(L^0)$ denotes the set of Mather (minimizing) measures of the Lagrangian $L^0$, where $L^0(x,v)=L(x,v,0)$. 

\begin{Lem}\label{new_gen_1}
For every Mather measure $\tmu\in\wtM(L^0)$, we have
\begin{align}\label{frac_00}
	\liminf_{\lambda\to 0^+}\int_{TM} u_\lambda(x)\frac{\partial L}{\partial  u }(x,v,0) \, d\tmu(x,v)\geq \int_{TM} V(x,0) \, d\tmu(x,v).
\end{align}
\end{Lem}
\begin{proof}
For each Mather measure $\tmu\in\wtM(L^0)$, we know that
\begin{equation*}
	\int_{TM}  L^0(x,v)\, d\tilde\mu(x,v)=\int_{TM}  L(x,v,0)\, d\tilde\mu(x,v)=-c(H^0).
\end{equation*}
Since $\tmu$ is necessarily a closed measure, invoking Proposition \ref{Pro_clomeasgeq0} we obtain 
    \begin{equation*}
		\int_{TM} L(x,v,\lambda u_\lambda(x))-\lambda V(x,\lambda)+c(H^0)\,d\tilde\mu(x,v)\geq 0,
    \end{equation*}
where $\lambda\in (0,\lambda_0]$. By combining the above results, we obtain
\begin{equation}\label{ineqintL}
\begin{aligned}
	0\leq \int_{TM}  L(x,v,\lam u_{\lam}(x))-\lam V(x,\lam)-L(x,v,0)\, d\tilde\mu(x,v). 
\end{aligned}  
\end{equation}
Let $S:=\textup{supp}(\tmu)$ denote the support of the measure $\tmu$, which is obviously a compact subset of $TM$. By property \eqref{d23rL4},
\[\left|L(x,v,\lam u_{\lam}(x))-L(x,v,0)-\lam u_\lam (x)\frac{\partial L}{\partial u }(x,v,0)\right|\leq \lambda\|u_{\lambda}\|_\infty \eta_{_S}(\lambda\|u_{\lambda}\|_\infty)\]
for all $(x,v)\in S$. Thus, \eqref{ineqintL} reduces to
\begin{equation*}
	0\leq \int_{TM} \lam u_\lam (x)\frac{\partial L}{\partial u }(x,v,0)+\lambda\|u_{\lambda}\|_\infty \eta_{_S}(\lambda\|u_{\lambda}\|_\infty)-\lam V(x,\lam)\, d\tmu(x,v).  
\end{equation*}
Dividing this inequality by $\lam$ and noting that $\eta_{_S}(\lambda\|u_{\lambda}\|_\infty)\to 0$ as $\lam\to 0^+$, we obtain 
$$
	0\leq \liminf_{\lam\to 0^+}\int_{TM}u_\lam(x)\frac{\partial L}{\partial u }(x,v,0) -V(x,\lam)\,d\tmu(x,v).
$$
Since $V(x,\lam)$ converges uniformly to $V(x,0)$, it thus verifies the desired inequality \eqref{frac_00}.
\end{proof}

As a corollary of Lemma \ref{new_gen_1}, we get a constraint on the possible accumulation points of $(u_\lambda)_{\lambda\in (0,\lambda_0]}$.
\begin{Cor}\label{Cor_constronlimit}
If for some sequence  $\lambda_k\to 0^+$, the viscosity solutions $u_{\lambda_k}$ converge uniformly to some function $u$, then 
\begin{align*}
\int_{TM} u(x)\frac{\partial L}{\partial  u }(x,v,0) \, d\tmu(x,v)\geq \int_{TM} V(x,0) \, d\tmu(x,v)
\end{align*}
for all Mather measures $\tmu\in\wtM(L^0)$.	
\end{Cor}

Next, consider the family $\mathcal{S}_0$ of continuous viscosity subsolutions $w:M\to\R$ of the critical equation \eqref{eq_G} that satisfies the inequality
\begin{align*}
	\int_{TM} \left[w(x)\frac{\partial L}{\partial  u }(x,v,0)-V(x,0)\right] \, d\tmu(x,v)\geq 0,\quad \text{for all $\tmu\in\wtM(L^0)$.}
\end{align*}
We introduce a continuous function $u_0:M\to\R$ defined by
\begin{equation}\label{u0defin}
	u_0(x):=\sup_{w\in \mathcal{S}_0} w(x).
\end{equation}
Clearly, $u_0$ itself is a viscosity subsolution of \eqref{eq_G} since it is the supremum of a family of viscosity subsolutions.	
\begin{Rem}
Corollary \ref{Cor_constronlimit} implies that every possible accumulation point $u$ of $(u_\lambda)_{\lambda\in (0,\lambda_0]}$ must satisfy  $u\leq u_0$. In the remainder of this section, we demonstrate that $u\geq u_0$, establishing the equality $u=u_0$ for any such accumulation point.
\end{Rem}

The following lemma will be very crucial in studying the relation between $u_\lambda$ and $u_0$. It involves the probability measure $\tmulax$ as defined in Definition \ref{def_measoncalicurve}.

\begin{Lem}\label{new_gen_2}
Consider any continuous viscosity subsolution $w:M\to\R$ of the critical equation \eqref{eq_G}. Then, for every $\lam\in (0,\lam_0]$ and $x\in M$, 
\begin{align}\label{impineq}
u_\lam(x)	\geq & w(x)-\frac{\int_{TM} \left[ w(y)\frac{\partial L}{\partial  u }(y,v,0)-V(y,\lam)+\lam^{-1}R_\lam (y,v)\right]\,d\tmulax(y,v)}{\int_{TM} \frac{\partial L}{\partial  u }(y,v,0)\,d\tmulax(y,v)}
\end{align}
where  $\tmulax$ is the probability measure on $TM$ defined in \eqref {def_pro_meas} and the term $R_\lam$ is
\begin{equation*}
	R_\lam(y,v)=L(y,v,\lam u_\lam(y))-L(y,v,0)-\lam u_\lam(y)\frac{\partial L}{\partial u}(y,v,0).
\end{equation*}
\end{Lem}

\begin{proof}
Let $w$ be a continuous viscosity subsolution of \eqref{eq_G}. Using a well-known approximation argument (see e.g. \cite{Fathi2012}), for any $\delta>0$, we can choose a smooth function $w_\delta \in C^\infty(M)$	 such that $\|w_\delta-w\|_\infty\leq \delta$ and
\begin{align*}
	H(y, d_yw_\delta,0)\leq c(H^0)+\delta,\quad \text{for all $y\in M$}.
\end{align*}
By the Fenchel inequality we have
\begin{equation}\label{d1hdf92}
	L(y,v,0)+c(H^0)\geq L(y,v,0)+H(y, d_yw_\delta,0)-\delta\geq   d_yw_\delta(v)-\delta, ~ \text{for all $(y,v)\in TM$}.
\end{equation}
Let $\xilax:(-\infty,0]\to M$ be the $u_\lam$-calibrated curve satisfying $\xilax(0)=x$. By Corollary \ref{Cor_differulam}, 
\begin{equation}\label{nue72j}
	\frac{d}{dt}u_\lambda(\xilax(t))=L\big(\xilax(t),\dxilax(t), \lam u_\lam(\xilax(t))\big)-\lam V(\xilax(t),\lam)+c(H^0),\quad\text{a.e. $t<0$},
\end{equation}
and by \eqref{d1hdf92} we get
\begin{equation}\label{fni10t}
\begin{aligned}
		\frac{d}{d t} w_\delta(\xilax(t))=&d_{\xilax(t)}w_\delta\,\big(\dxilax(t)\big)
	\leq  L(\xilax(t),\dxilax(t),0)+c(H^0)+\delta, \quad\text{a.e. $t<0$}.
\end{aligned}
\end{equation}
Subtracting  \eqref{fni10t} from  \eqref{nue72j} yields
\begin{align*}
	\frac{d}{d t}u_\lam(\xilax(t))
	\geq & \frac{d}{d t} w_\delta(\xilax(t)) +L\big(\xilax(t),\dxilax(t), \lam u_\lam(\xilax(t))\big)-L(\xilax(t),\dxilax(t),0)\\
	&\qquad\qquad\qquad \qquad\qquad \qquad\qquad \qquad\qquad\qquad\qquad-\lam V(\xilax(t),\lam)-\delta
\end{align*}
a.e. $t<0$. This thus gives
\begin{equation}\label{ncqo2}
\begin{aligned}
	 \frac{d}{d t}u_\lam(\xilax(t))-\lam u_\lam (\xilax(t)) & \frac{\partial L}{\partial  u }(\xilax(t),\dxilax(t),0)\\
	 &\geq \frac{d}{d t} w_\delta(\xilax(t))+R_\lam(\xilax(t),\dxilax(t))-\lam V(\xilax(t),\lam)-\delta,
 \end{aligned}
 \end{equation}
a.e. $t<0$, where the continuous function $R_\lam(\cdot,\cdot): TM\to \R$ is given by
\begin{equation*}
	R_\lam(y,v)=L(y,v,\lam u_\lam(y))-L(y, v,0)-\lam  u_\lam(y) \frac{\partial L}{\partial u}(y,v,0).
\end{equation*}

Multiplying both sides of inequality \eqref{ncqo2} by an integrating factor $e^{\lambda\int^0_t \frac{\partial L}{\partial u }(\xilax(s),\dxilax(s),0)\, ds}$, we integrate over $(-T,0]$ and use integration by parts to obtain
\begin{align*}
	u_\lam(x)- u_\lam(\xilax(-T)) & e^{\lam\int^0_{-T} \frac{\partial L}{\partial u }(\xilax(s),\dxilax(s),0)\, ds}\\
	 \geq & w_\delta(x)-w_\delta (\xilax(-T)) e^{\lam\int^0_{-T} \frac{\partial L}{\partial u }(\xilax(s),\dxilax(s),0)\, ds}\\
	&-\int_{-T}^0   w_\delta(\xilax(t)) \frac{d}{d t}\left(e^{\lam\int^0_t \frac{\partial L}{\partial u }(\xilax(s),\dxilax(s),0)\, ds}\right)dt\\
	&+ \int_{-T}^0 \left(R_\lam(\xilax(t),\dxilax(t))-\lam V(\xilax(t),\lam)-\delta\right)e^{\lam\int^0_t \frac{\partial L}{\partial u }(\xilax(s),\dxilax(s),0)\, ds}\,dt.	
\end{align*}
Sending $T\to +\infty$ and using Lemma \ref{dw72je}, we get  
\begin{equation}\label{gdht}
u_\lam(x)	\geq   w_\delta(x) +\mathbf{E}_1+\mathbf{E}_2
\end{equation} 
where the terms $\mathbf{E}_1$ and $\mathbf{E}_2$ are given by
\begin{align*}
 &\mathbf{E}_1=-\int_{-\infty}^0   w_\delta(\xilax(t)) \frac{d}{d t}\left(e^{\lam\int^0_t \frac{\partial L}{\partial u }(\xilax(s),\dxilax(s),0)\, ds}\right) dt	,\\
 &\mathbf{E}_2=\int_{-\infty}^0 \left(R_\lam(\xilax(t),\dxilax(t))-\lam V(\xilax(t),\lam)-\delta\right)e^{\lam\int^0_t \frac{\partial L}{\partial u }(\xilax(s),\dxilax(s),0)\, ds}\,dt.
\end{align*}
For the term $\mathbf{E}_1$, using \eqref{useident0} and the probability measure $\tmulax$ defined in \eqref {def_pro_meas}, we find
\begin{align*}
	\mathbf{E}_1=&
	\lam \int_{-\infty}^0   w_\delta(\xilax(t))\, \frac{\partial L}{\partial  u }(\xilax(t),\dxilax(t),0)\,  e^{\lam\int^0_t \frac{\partial L}{\partial u }(\xilax(s),\dxilax(s),0)\, ds}\,dt\\
	=&\lam \left(\int_{-\infty}^0   e^{\lam\int^0_t  \frac{\partial L}{\partial  u }(\xilax(s),\dxilax(s),0)\, ds}\, dt\right) \int_{TM}   w_\delta(y) \frac{\partial L}{\partial u }(y,v,0)\,d\tmulax(y,v).
\end{align*}
Invoking Lemma \ref{Lem_measidentity} to this equality gives
\begin{equation}\label{jaydw21}
	\mathbf{E}_1=\frac{-1}{\int_{TM} \frac{\partial L}{\partial  u }(y,v,0)\,d\tilde{\mu}_x^\lam(y,v)}\int_{TM}   w_\delta(y) \,\frac{\partial L}{\partial u }(y,v,0)\,d\tmulax(y,v).
\end{equation}
For the term $\mathbf{E}_2$, using the probability measure $\tmulax$ and Lemma \ref{Lem_measidentity} again we obtain 
\begin{equation}\label{jaydw22}
\begin{aligned}
\mathbf{E}_2
	=&  \left(\int_{-\infty}^0   e^{\lam\int^0_t  \frac{\partial L}{\partial  u }(\xilax(s),\dxilax(s),0)\, ds}\, dt\right)  \int_{TM}  R_\lam(y,v)-\lam V(y,\lam)-\delta\,d\tmulax(y,v)\\
	=&\frac{-1}{\int_{TM} \frac{\partial L}{\partial  u }(y,v,0)\,d\tilde{\mu}_x^\lam(y,v)}\int_{TM} \lam^{-1} R_\lam (y,v)- V(y,\lam)-\lam^{-1}\delta\,d\tmulax(y,v).
\end{aligned}	
\end{equation}

Finally, combining the expressions \eqref{gdht} and \eqref{jaydw21}-\eqref{jaydw22},  we get
\[
u_\lam(x)	\geq  w_\delta(x)-\frac{\int_{TM}  w_\delta(y) \,\frac{\partial L}{\partial  u }(y,v,0)+\lam^{-1}R_\lam(y,v)-V(y,\lam)-\lam^{-1}\delta\,d\tmulax(y,v)}{\int_{TM} \frac{\partial L}{\partial  u }(y,v,0)\,d\tmulax(y,v)}.
\]
Sending $\delta\to 0$ in this estimate, we obtain the desired inequality, proving the lemma.
\end{proof}

\subsection{Proof of the convergence result}
We are now ready to complete the proof of Theorem \ref{mainresult1}.

\begin{proof}[Proof of Theorem \ref{mainresult1}]
Let us recall the continuous function $u_0: M\to \R$ given in 
\eqref{u0defin}:
\begin{equation*}
	u_0(x)=\sup_{w\in\mathcal{S}_0}w(x)
\end{equation*}
where $\mathcal{S}_0$ is the family of continuous viscosity subsolutions $w$ of  equation \eqref{eq_G} satisfying 
\begin{align}\label{ccsub}
	\int_{TM} \left[w(x)\frac{\partial L}{\partial  u }(x,v,0)-V(x,0)\right] \, d\tilde\nu(x,v)\geq 0,\quad \text{for all $\tilde\nu\in\wtM(L^0)$.}
\end{align}
We aim to show that $u_0$ is the unique limit of the family $(u_\lambda)_{\lambda\in(0,\lambda_0]}$. This is equivalent to proving that any convergent subsequence has $u_0$ as its limit. Consider a subsequence $\lambda_k\to 0$ for which  $u_{\lam_k}$ converges uniformly to some function $u$. According to Corollary \ref{Cor_constronlimit}, this function $u$ belongs to $\mathcal{S}_0$, implying $u(x)\leq u_0(x)$ for all $x\in M$. 

To establish the reverse inequality, we make use of Lemma \ref{new_gen_2}. Indeed, let   $x\in M$ and $w\in \mathcal{S}_0$, it follows from Lemma \ref{new_gen_2} that 
 \begin{align}\label{ajeiq1}
u_{\lam_k}(x)	\geq & w(x)-\frac{\int_{TM}  \left[w(y) \,\frac{\partial L}{\partial  u }(y,v,0)-V(y,\lam_k)+\lam_k^{-1}R_{\lambda_k}(y,v)\right]\,d\tilde{\mu}_x^{\lam_k}(y,v)}{\int_{TM} \frac{\partial L}{\partial  u }(y,v,0)\,d\tilde{\mu}_x^{\lam_k}(y,v)}
\end{align}
where
\begin{equation*}
	R_{\lam_k}(y,v)=L(y,v,\lam_k u_{\lam_k}(y))-L(y,v,0)-\lam_k u_{\lam_k}(y)\frac{\partial L}{\partial u}(y,v,0).
\end{equation*}
Recall that by Proposition \ref{Pro_weakconvergence} the probability measures $\tilde{\mu}^{\lam_k}_x$ have support contained in a common compact subset $S$ of $TM$, so using property \eqref{d23rL4} we get 
\begin{align}
	\left|\int_{TM} \lam_k^{-1}R_{\lambda_k}(y,v)\,d\tilde{\mu}_x^{\lam_k}(y,v)\right| &\leq  \left|   \int_{TM} \|u_{\lam_k}\|_\infty\,\eta_{_S}(\lam_k\|u_{\lam_k}\|_\infty)\,d\tilde{\mu}_x^{\lam_k}(y,v)   \right|\nonumber\\
	& \leq \|u_{\lam_k}\|_\infty\,\eta_{_S}(\lam_k\|u_{\lam_k}\|_\infty)\nonumber\\
	& \rightarrow 0 \quad\text{as $\lam_k\to 0^+$.}\label{ajeiq2}
\end{align}
Meanwhile, from Proposition \ref{Pro_weakconvergence} we can assume (extracting a further subsequence if necessary) that $\tilde{\mu}^{\lam_k}_x$  converges weakly  to a Mather measure $\tilde\mu$. Consequently, sending $\lambda_k\to 0^+$ in \eqref{ajeiq1}, combined with \eqref{ajeiq2} and the fact that $V(\cdot,\lam_k)$ converges to $V(\cdot,0)$ on $M$, we get
\begin{equation}\label{dnao1}
	u(x)\geq w(x)-\frac{\int_{TM}  w(y) \,\frac{\partial L}{\partial  u }(y,v,0)-V(y,0)\,d\tilde{\mu}(y,v)}{\int_{TM} \frac{\partial L}{\partial  u }(y,v,0)\,d\tilde{\mu}(y,v)}.
\end{equation}
In view of \eqref{ccsub} and $\frac{\partial L}{\partial u}(y,v,0)<0$, inequality \eqref{dnao1} directly leads to $u(x)\geq w(x)$. Since $w\in \mathcal{S}_0$ is arbitrary, it follows that 
\[u(x)\geq \sup_{w\in\mathcal{S}_0}w(x)=u_0(x).\]
Combining this with the earlier inequality $u(x)\leq u_0(x)$, we conclude that  $u(x)=u_0(x)$.
This finally completes the proof.	
\end{proof}

%%%%%%%%%%%%%%%%%%%%%%%%%%%%%%%%%%%%%%%%%%%%%%%%%%%%%%%%%%%%%%%%%%

\section{Another formula for the limit solution}\label{section_Anotherformula}
In this section, we provide an alternative characterization of the limit solution by proving Theorem \ref{mainresult2}. This theorem expresses $u_0$ using the Peierls barrier, Mather measures, and the potential $V(x,0)$. The significance of this formula lies in its flexibility: by varying the potential $V$, we can select a wider variety of solutions to the critical Hamilton-Jacobi equation.

Let $h(\cdot,\cdot): M\times M\to\R $ denote the Peierls barrier of the Lagrangian $L^0(x,v)$. We begin with the following lemma.

\begin{Lem}\label{lem_anoform1}
	For every $x\in M$, the limit solution $u_0(x)$ satisfies:
	\begin{equation*}
		u_0(x)\leq \frac{\int_{TM}  h(y,x)\frac{\partial L}{\partial u}(y,v,0)+V(y,0)\,d\tilde\mu(y,v)}{\int_{TM} \frac{\partial L}{\partial u}(y,v,0) \,d\tmu(y,v)},\
	\end{equation*}
	for any Mather measure $\tmu\in \wtM(L^0)$.
\end{Lem}

\begin{proof}
	Let $x\in M$. Since $u_0$ is a viscosity solution of \eqref{eq_G}, we can apply Proposition \ref{prop h} to obtain that for any $y\in M$,  $u_0(x)\leq u_0(y)+h(y,x)$. Multiplying both sides of this inequality by $\frac{\partial L}{\partial u}(y,v,0)<0$ gives
\begin{equation*}
	u_0(x) \frac{\partial L}{\partial u}(y,v,0) \geq u_0(y)\frac{\partial L}{\partial u}(y,v,0)+ h(y,x)\frac{\partial L}{\partial u}(y,v,0),\quad \text{for all $(y,v)\in TM$.}
\end{equation*}
Then, for each Mather measure $\tilde\mu \in \wtM(L^0)$, it follows that
\begin{align*}
	u_0(x) \int_{TM} \frac{\partial L}{\partial u}(y,v,0)\,d\tilde\mu(y,v)
	& \geq \int_{TM} u_0(y)\frac{\partial L}{\partial u}(y,v,0)+ h(y,x)\frac{\partial L}{\partial u}(y,v,0)\,d\tilde\mu(y,v).
\end{align*}
Using Corollary \ref{Cor_constronlimit}, the limit $u_0$ satisfies
\begin{align*}
\int_{TM} u_0(y)\frac{\partial L}{\partial  u }(y,v,0) \, d\tmu(y,v)\geq \int_{TM} V(y,0) \, d\tmu(y,v).
\end{align*} 
Combing these inequalities gives
\begin{align}\label{p0wq8}
	u_0(x) \int_{TM} \frac{\partial L}{\partial u}(y,v,0)\,d\tilde\mu(y,v)
     \geq \int_{TM}  h(y,x)\frac{\partial L}{\partial u}(y,v,0)+V(y,0)\,d\tilde\mu(y,v).
\end{align}
Dividing both sides of \eqref{p0wq8} by $ \int_{TM} \frac{\partial L}{\partial u}(y,v,0)\,d\tilde\mu(y,v)$ (which is negative), we obtain the asserted inequality.
\end{proof}

Inspired by Lemma \ref{lem_anoform1}, we define a function $\Theta:M\to\R$ by 
	\begin{equation*}
		\Theta(x):=\inf_{\tilde{\mu}\in \wtM(L^0)}  \frac{\int_{TM}  h(y,x)\frac{\partial L}{\partial u}(y,v,0)+V(y,0)\,d\tilde\mu(y,v)}{\int_{TM} \frac{\partial L}{\partial u}(y,v,0) \,d\tilde\mu(y,v)}.
    \end{equation*}

\begin{Lem}\label{lem_0sai21}
The function $\Theta$ is a continuous viscosity subsolution of equation \eqref{eq_G}.
\end{Lem} 
\begin{proof}
For each Mather measure $\tilde\mu\in \wtM(L^0)$, we define a function $h_{\tilde\mu}: M\to \R$ by
\[h_{\tilde\mu}(x):= \frac{\int_{TM}  h(y,x)\frac{\partial L}{\partial u}(y,v,0)\, d\tilde\mu(y,v)}{\int_{TM} \frac{\partial L}{\partial u}(y,v,0) \,d\tilde\mu(y,v)   }.\]
Then, $h_{\tilde\mu}$ is a continuous viscosity subsolution of equation \eqref{eq_G},  because it is a convex combination of equi-Lipschitz viscosity solutions $\{h_y(x):=h(y,x)\}_{y\in M}$ of \eqref{eq_G}. Next, define a function $h_{\tilde\mu, V}: M\to \R$ by
\[
h_{\tilde\mu, V}(x):=h_{\tilde\mu}(x)+\frac{\int_{TM} V(y,0)\,d\tilde\mu(y,v)}{\int_{TM} \frac{\partial L}{\partial u}(y,v,0) \,d\tilde\mu(y,v)}.
\]
Clearly, $h_{\tilde\mu, V}$ is also a viscosity subsolution of \eqref{eq_G} since it differs from $h_{\tilde\mu}$ only by a constant. Therefore, as the infimum of a family of equi-Lipschitz subsolutions,  $\Theta=\inf_{\tilde{\mu}\in \wtM(L^0)} h_{\tilde\mu, V}$ itself is a continuous viscosity subsolution of \eqref{eq_G}.
\end{proof}

Our next statement describes that $\Theta(x)$ is indeed the limit solution $u_0(x)$, which proves Theorem \ref{mainresult2}.

\begin{The}
	The limit solution $u_0:M\to\R$ obtained in Theorem \ref{mainresult1} is equal to $\Theta$.
\end{The} 
\begin{proof}
We first note that by Lemma \ref{lem_anoform1}, $u_0(x)\leq \Theta(x)$ for all $x\in M$. It remains to check the opposite inequality $u_0\geq \Theta$.

Fix $x\in M$ and define a function $w:M\to \R$ by 
\[w(y):=-h(y,x).\] 
By Proposition \ref{prop h}, $w$ is actually a viscosity subsolution of equation \eqref{eq_G}. Applying Lemma \ref{new_gen_2} to the function $w$, we obtain
\begin{align*}
	u_\lambda(x)\geq-h(x,x)-\frac{\int_{TM}  \left[-h(y,x)\frac{\partial L}{\partial u}(y,v,0) -V(y,\lam)+\lam^{-1}R_\lam (y,v)\right]\,d\tilde{\mu}_x^\lam(y,v)}{\int_{TM} \frac{\partial L}{\partial u}(y,v,0) \,d\tilde{\mu}_x^\lam(y,v) }.
\end{align*}
Sending $\lambda$ to zero and using argument similar to the proof of \eqref{dnao1}, we get 
\begin{align}\label{mw3yw}
	u_0(x)\geq -h(x,x)+\frac{\int_{TM}  h(y,x)\frac{\partial L}{\partial u}(y,v,0)+ V(y,0)\,d\tilde\nu(y,v)}{\int_{TM} \frac{\partial L}{\partial u}(y,v,0) \,d\tilde\nu(y,v)   }.
\end{align}
for some Mather measure $\tilde\nu$. 

In particular, for $x$ in the projected Aubry set $\cA_{L^0}$ of $L^0$,  we have $h(x,x)=0$, and hence \eqref{mw3yw} directly yields 
\begin{align*}
	u_0(x)\geq \frac{\int_{TM}  h(y,x)\frac{\partial L}{\partial u}(y,v,0)+ V(y,0)\,d\tilde\nu(y,v)}{\int_{TM} \frac{\partial L}{\partial u}(y,v,0) \,d\tilde\nu(y,v) }
	\geq \Theta(x),\quad \text{for $x\in \cA_{L^0}$}.
\end{align*}
Since $\Theta$ is a viscosity subsolution (see Lemma \ref{lem_0sai21}) of equation \eqref{eq_G}, we conclude from Proposition \ref{Aubry_uniqueness} that $u_0\geq \Theta $ everywhere on $M$. This completes the proof.
\end{proof}

We conclude this section by giving a straightforward example that explicitly illustrates how the limit solution $u_0$ depends on $V$.
\begin{Ex}\label{ex_simexamp2}
	Consider a Hamiltonian $H(x,p,u)=u+\frac{1}{2}\|p+\alpha\|^2$ on $T^*\T^n\times\R$, where $\T^n=\R^n/\Z^n$ and $\alpha\in \R^n$ is a Diophantine vector. Note that the critical value $c(H^0)=\frac{1}{2}\|\alpha\|^2$, so the corresponding perturbation problem is given by
	\begin{equation}\label{neq1ie}
		\lambda u_\lambda(x)+\frac{1}{2}\|D u_\lambda+\alpha\|^2+\lambda V(x,\lambda)=\frac{1}{2}\|\alpha\|^2,
	\end{equation}
	which approximates the critical H-J equation
	\begin{equation}\label{neq2ie}
		\frac{1}{2}\|Du+\alpha\|^2=\frac{1}{2}\|\alpha\|^2.
	\end{equation}
\begin{itemize}
	\item In such a case, the projected Mather set for the critical Hamiltonian $H^0(x,p)$ is $\mathcal M=\T^n$, on which the dynamics is a Diophantine linear flow. The corresponding Peierls barrier $h(x,y)\equiv 0$ on $M\times M$. Consequently, all solutions of the critical equation \eqref{neq2ie} are {\em constant}.
	\item On the other hand, the projected Mather measure is unique and is the Lebesgue measure on $\mathcal M$, so Theorem \ref{mainresult2} yields that the limit solution $u_0$ is of the form 
\begin{equation}
	u_0(x)\equiv\int_{\T^n} V(\theta,0)\,d\theta,\quad \forall~x.
\end{equation}
Thus, by choosing suitable potential $V$, we can generate any desired solution to equation \eqref{neq2ie}. 
\end{itemize}
\end{Ex}

%%%%%%%%%%%%%%%%%%%%%%%%%%%%%%%%%%%%%%%%%%%%%%%%%%%%%%%%%%%%%%%%%%
\section{Applications of the new selection principle}\label{section_proofofCandD}
In this section, we prove Theorem \ref{mainresult3}, Theorem \ref{mainresult4} and Theorem \ref{mainresult6}, highlighting the versatility of the new selection principle in addressing issues such as the variational characterization of viscosity solutions, the refinement of uniqueness results, and the structural properties of equilibrium Mather measures.

\subsection{Proof of Theorems \ref{mainresult3} and \ref{mainresult4}}

Theorem \ref{mainresult3} reveals the properties of the operator $\Psma$ (see Definition \ref{def_operaP}) and its role in characterizing solutions to the  Hamilton-Jacobi equation \[G(x, d_xu)=c(G) \quad \text{in $M$}.\]
Its proof relies on applying Theorems \ref{mainresult1}--\ref{mainresult2} to the following perturbation problem:
\begin{equation*}
	\sigma(x)\lambda u_\lambda(x)+G(x, d_xu_\lambda)-\lambda V(x,\lambda)=c(G)\quad \text{in $M$.}
\end{equation*}

To prove  Theorem  \ref{mainresult3}, we recall the definition of the operator $\Psma$ for the Hamiltonian $G$: 
	given a continuous function $\sigma: M\to (0,+\infty)$, we set \[\Psma: C(M)\to  C(M)\] 
	as follows: for every $\varphi\in C(M)$, the function $\Psma\varphi:M\to \R$ is
\begin{equation*}
\begin{aligned}
	\Psma\varphi(x):=&\inf_{\tmu \in \tilde{\mathfrak{M}}(L_G)} \left\{ \frac{\int_{TM}  \sigma(y) h_G(y,x)\,d\tmu(y,v)}{\int_{TM} \sigma(y) \,d\tmu(y,v)   } +\frac{\int_{TM}  \sigma(y)\varphi(y)\,d\tmu(y,v)}{\int_{TM} \sigma(y) \,d\tmu(y,v)} \right\}\\
	=&\inf_{\mu\in \mathfrak{M}(L_G)} \left\{ \frac{\int_{M}  \sigma(y) h_G(y,x)\,d\mu(y)}{\int_{M} \sigma(y) \,d\mu(y)   } +\frac{\int_M  \sigma(y)\varphi(y)\,d\mu(y)}{\int_{M} \sigma(y) \,d\mu(y) } \right\}
\end{aligned}
    \end{equation*}
where $ \tilde{\mathfrak{M}}(L_G)$ (resp. $ \mathfrak{M}(L_G)$) is the set of all  (projected) Mather measures of $L_G$.

\begin{proof}[Proof of Theorem \ref{mainresult3}] \textbf{Assertion (i): Lipschitz continuity of $\Psma$.} 

We first notice that by the definition of $\Psma$, 
\[\Psma(\varphi+c)=\Psma(\varphi)+c,\quad\forall~ \varphi\in C(M)~\text{and}~ c\in \R .\]
Moreover, it is easy to find that $\Psma$ is order preserving: $\Psma f\leq \Psma g$ whenever 
$f\leq g$. Consider two functions $\varphi_1, \varphi_2\in C(M)$, we have
\[ \varphi_2-\|\varphi_1-\varphi_2\|_\infty  \leq \varphi_1\leq \varphi_2+\|\varphi_1-\varphi_2\|_\infty.\]
Using the order preserving property of $\Psma$, we obtain
\[    \Psma\varphi_2-\|\varphi_1-\varphi_2\|_\infty  \leq \Psma\varphi_1\leq \Psma\varphi_2+\|\varphi_1-\varphi_2\|_\infty.\]
This yields assertion (i).

\textbf{Assertion (ii): Solutions to the equation $G(x, d_xu)=c(G)$.} 

The proof relies on a specific perturbation of the equation $G(x, d_xu)=c(G)$ to which we can apply Theorem \ref{mainresult1} and Theorem \ref{mainresult2}. Indeed, we study the following perturbed equation for $\lambda>0$:
\begin{equation}\label{akjdiuqwi}
	\sigma(x)\lambda u_\lambda(x)+G(x, d_xu_\lambda)-\lambda V(x,\lambda)=c(G) \quad \text{in $M$}
\end{equation}
where $V(x,\lambda)$ converges uniformly to $V(x,0)$ as $\lambda\to 0$, with $V(x,0)$ satisfying
\[V(x,0)=-\sigma(x)\varphi(x), \quad\forall~x\in M.\] 
This corresponds to the equation \eqref{eq_H_lambda} by setting $H(x,p,u)=\sigma(x) u+ G(x,p)$. Such a Hamiltonian $H$ satisfies conditions {\bf(H1)-(H4)}, so we can      
apply Theorem \ref{mainresult1} to obtain that the solution $u_\lambda$ of equation \eqref{akjdiuqwi} converges uniformly to a limit function $u_0:M\to \R$ as $\lambda\to 0$. Moreover, $u_0$ is a viscosity solution of 
\begin{equation*}
G(x, d_xu_\lambda)=c(G)\quad \text{in $M$}.
\end{equation*}
Theorem \ref{mainresult2} further gives
\begin{align*}
	u_0(x)= &\inf_{\mu\in \mathfrak{M}(L_G)}  \frac{\int_{TM}  -\sigma(y)h_G(y,x)\,d\mu(y)+\int_M V(y,0)\,d\mu(y)}{\int_{TM} -\sigma(y) \,d\mu(y)   }\\
	=&\inf_{\mu\in \mathfrak{M}(L_G)} \left\{ \frac{\int_{M}  \sigma(y) h_G(y,x)\,d\mu(y)}{\int_{M} \sigma(y) \,d\mu(y)   } +\frac{\int_M  \sigma(y)\varphi(y)\,d\mu(y)}{\int_{M} \sigma(y) \,d\mu(y) } \right\}
\end{align*}
where $\mathfrak{M}(L_G)$ is the set of all projected Mather measures of $L_G$.
The last line of the above equation is precisely $\Psma\varphi(x)$. This confirms that $\Psma\varphi$ is indeed a viscosity solution of $G(x, d_xu_\lambda)=c(G)$, thus proving assertion (ii).

\textbf{Assertion (iii): Fixed Point Characterization.} 

We now prove that a continuous function $u:M\to\R$ is a viscosity solution of the  equation $G(x,d_xu)=c(G)$ if and only if it is a fixed point of the operator $\Psma$.

Suppose that $u=\Psma u$. From assertion (ii), it follows that $u$ is a viscosity solution of the equation $G(x,d_xu)=c(G)$.

Conversely, suppose that $u:M\to \R$ is a viscosity solution of the H-J equation $G(x,d_xu)=c(G)$. We apply Theorem \ref{mainresult1} and Theorem \ref{mainresult2} to the following approximation of the equation $G(x,d_xu)=c(G)$:
\begin{equation}\label{daqiq12d}
	\lam \sigma(x) u_\lam(x)+G(x,d_x u_\lam)+\lam U(x,\lambda) =c(G) \quad \text{in $M$}
\end{equation}
where $U(x,\lambda)$ is chosen so that $U(x,0)=-\sigma(x)u(x)$ and 
\begin{equation}\label{caq7h}
	 \sup_{x\in M}|U(x,\lambda)-U(x,0)|\leq \lambda. 
\end{equation}
Notably, $U(x, \lambda)$ converges uniformly to $U(x,0)$. Then, analogous to the proof of assertion (ii), we infer that the solution $u_\lambda$ of \eqref{daqiq12d} satisfies 
\begin{equation}\label{d18je}
	 \lim_{\lambda\to 0} u_\lambda(x)=\Psma u(x),\quad \forall~x\in M.
\end{equation}
On the other hand,  consider two functions on $M$ defined by 
\[u_\lambda^-(x):=u(x)-\frac{\lambda}{\min_{x\in M} \sigma(x)},\qquad u_\lambda^+(x):=u(x)+\frac{\lambda}{\max_{x\in M} \sigma(x)}.\]
Due to \eqref{caq7h}, $u_\lambda^-$ and $u_\lambda^+$ are respectively viscosity subsolution and supersolution of equation \eqref{daqiq12d}. By the comparison principle, we get
\begin{equation*}
	u_\lambda^-(x)\leq u_\lambda(x)\leq u_\lambda^+(x),
\end{equation*}
from which we infer that 
\begin{equation*}
	u(x)=\lim_{\lambda\to 0} u^-_\lambda(x)\leq\lim_{\lambda\to 0} u_\lambda(x)\leq\lim_{\lambda\to 0} u^+_\lambda(x)=u(x),\quad \forall~x\in M.
\end{equation*}
This combined with \eqref{d18je} verifies that $u=\Psma u$.
\end{proof}

By Theorem \ref{mainresult3}, the image of the operator
	$\Psma$ in $C(M)$ equals the set of all viscosity solutions to $G(x, d_xu)=c(G)$. 
Thus, our selection principle (combined effects of both the vanishing discount process and potential perturbations) enables to select any desired viscosity solution to $G(x, d_xu)=c(G)$. Meanwhile, Theorem \ref{mainresult3} also leads naturally to the proof of Theorem \ref{mainresult4}.

\begin{proof}[Proof of Theorem \ref{mainresult4}] 
Let $u_1, u_2\in C(M)$ be two viscosity solutions of the H-J equation $G(x,d_xu)=c(G)$. We have seen from Theorem \ref{mainresult3} that $u_1, u_2$ are both fixed points of the operator $\Psma$, so for each $x\in M$,
\begin{equation*}
	u_i(x):=\inf_{\mu\in \mathfrak{M}(L_G)} \left\{ \frac{\int_{M}  \sigma(y) h_G(y,x)\,d\mu(y)}{\int_{M} \sigma(y) \,d\mu(y)   } +\frac{\int_M  \sigma(y)u_i(y)\,d\mu(y)}{\int_{M} \sigma(y) \,d\mu(y) } \right\},\quad i=1,2.
\end{equation*}
Therefore, if \[
	\int_M \sigma(y) u_1(y)\,d\mu\leq \int_M \sigma(y) u_2(y)\,d\mu,\quad \text{for any}~\mu\in \mathfrak{M}(L_G),
\]
then it follows directly from the above representation that $u_1\leq u_2$ everywhere on $M$. This completes the proof of Theorem \ref{mainresult4}.
\end{proof}

\subsection{Proof of Theorem \ref{mainresult6}}

Finally, we turn to study  measures that realize the infimum in the definition of the operator $\Psmone:=\Psma$ with $\sigma\equiv 1$. In particular, we study the generic properties of the equilibrium Mather measures.

\begin{defn}
Given  $\varphi\in C(M)$ and $x\in M$. A Mather measure $\tilde\mu\in \tilde{\mathfrak{M}}(L_G)$ on $TM$ is called 
 an {\em equilibrium Mather measure} for $\varphi$ and $x$ if it satisfies:  
\[ \int_{TM}  h_G(y, x)\,d\tmu(y,v) +\int_{TM} \varphi(y)\,d\tmu(y,v)=\Psmone\varphi(x).\]
Hence, the projection $\mu=\pi_{\#}\tmu$ of $\tmu$ on $M$ satisfies
\[ \int_{M}  h_G(y, x)\,d\mu(y) +\int_{M} \varphi(y)\,d\mu(y)=\Psmone\varphi(x),\]
and we call $\mu$ the {\em equilibrium projected Mather measure} for $\varphi$ and $x$. We denote by $\tilde{\mathcal{E}_{G}}(\varphi, x)$ the set of all equilibrium Mather measures for $\varphi$ and $x$, and by  $\mathcal{E}_{G}(\varphi, x)$  the set of all equilibrium projected Mather measures for $\varphi$ and $x$. 
\end{defn}

\begin{proof}[Proof of Theorem \ref{mainresult6}]
\textbf{(i)} The Lagrangian $L_G$ associated with $G$ is a $C^2$ Tonelli Lagrangian. Classical Aubry-Mather theory asserts that Mather measures are invariant under the Euler-Lagrange flow, and the set of all Mather measures $\tilde{\mathfrak{M}}(L_G)$ is a compact and convex set whose extremal points are ergodic invariant measures. It follows directly that $\tilde{\mathcal{E}_{G}}(\varphi, x)$, the set of all equilibrium Mather measures for $\varphi$ and $x$, is also convex and compact.

Let $\tilde\mu$ be an extremal point of the convex set $\tilde{\mathcal{E}_{G}}(\varphi, x)$. If $\tilde\mu$ is not ergodic, then we may write $\tilde\mu= k\tilde\mu_1+(1-k)\tilde\mu_2$, where  $k\in (0,1)$ and $\tilde\mu_1, \tilde\mu_2\in \tilde{\mathfrak{M}}(L_G)$ are ergodic Mather measures. Then it follows that 
\begin{align*}
	\Psmone\varphi(x)=&\int_{TM}  h_G(y, x) + \varphi(y)\,d\tmu(y,v)\\
	=&k\left[\int_{TM}  h_G(y, x) + \varphi(y)\,d\tmu_1(y,v)\right]+(1-k)\left[\int_{TM}  h_G(y, x) + \varphi(y)\,d\tmu_2(y,v)\right]\\
	\geq & k \Psmone\varphi(x)+(1-k)\Psmone\varphi(x)=\Psmone\varphi(x).
\end{align*}
This implies $\tilde\mu_1, \tilde\mu_2\in \tilde{\mathcal{E}_{G}}(\varphi, x)$, and it contradicts that $\tilde\mu$ is an extremal point of $\tilde{\mathcal{E}_{G}}(\varphi, x)$.

\textbf{(ii)} By Mather's graph property 
(Theorem \ref{appendthm01}), the linear map $\pi_\#: \wtM(L_G)\to \mathfrak{M}(L_G) $ is injective. This implies that $\operatorname{card}(\tilde{\mathcal{E}_{G}}(\varphi, x))=\operatorname{card}(\mathcal{E}_{G}(\varphi, x))$. Therefore,  it suffices to prove the generic uniqueness for $\mathcal{E}_{G}(\varphi, x)$. To begin with we shall work in the following setting .
	\begin{itemize}
  \item  Let $E=C(M)$ equipped with the supremum norm. $E$ is a Banach space.

  \item  Let $F$ be the space of of finite Borel measures on $M$ endowed with the weak-$*$ topology. In particular, for $\mu_k, \mu\in F$,
      $$\lim\limits_{k\to \infty}\,\mu_k=\mu \Longleftrightarrow \lim\limits_{k\to\infty}\int_{M}f\, d\mu_k=\int_{M}f d\mu,\quad\forall f\in E.$$
  Moreover,  $F$ is  metrizable and separable.
 
  \item For any fixed $x\in M$, define a linear map $\mathcal{L}_x: F\to \R$ by setting 
  \[\mathcal{L}_x(\mu):=\int_M h(y,x)\, d\mu(y), ~~\forall~\mu\in F.\]
  \item  Let $\Phi: E\to F^*$ be a linear map such that for $\varphi\in E$, $\Phi(\varphi)\in F^*$ and
    \[\langle\Phi(\varphi),\mu\rangle=\int_M \varphi\, d\mu,~ \mu\in F.\]
     \item  Let $K=\mathfrak{M}(L_G)$ be the set of all projected Mather measures for the Lagrangian $L_G$. It is a convex, compact subset of $F$. Observe that $K$ is separated by $E$, namely, for every pair $\mu_1, \mu_2\in K$, there exists $f\in E$ such that 
  \[   \int_M f\,d\mu_1\neq  \int_M f\,d\mu_2. \]
\end{itemize}
Fix a point $x\in M$ and a function $\varphi\in E$. We denote the value
\begin{align*}
	m(\varphi,x)= &\inf_{\mu\in K}\Big( \mathcal{L}_x(\mu)+\langle\Phi(\varphi),\mu\rangle  \Big)\\
	= & \inf_{\mu\in K} \int_{M} h(y,x)+\varphi(y)\,d\mu(y).
\end{align*} 
As a consequence, the considered set $\mathcal{E}_{G}(\varphi, x)$ is 
\[ \mathcal{E}_{G}(\varphi, x)=\{ ~\mu\in K ~|~ \mathcal{L}_x(\mu)+\langle\Phi(\varphi),\mu\rangle=m(\varphi,x) ~\}.  \]
Then, by the genericity approach of Ma\~n\'e in \cite[Section 3]{Mane1996}, see also \cite[Section 2]{Bernard_Contreras2008}, it follows that there exists a residual subset $\mathcal{O}_x\subset E=C(M)$ such that for any $\varphi\in \mathcal{O}_x$, the set $\mathcal{E}_{G}(\varphi, x)$ contains just one element. This proves the generic uniqueness result.

\textbf{(iii)} Consider a countable dense subset $\mathcal{D}=\{x_n\}_{n\geq 1}$ of the manifold $M$, by the above result we have a residual subset  $\mathcal{O}_{x_n}C(M)$ for each $n$. Thus, the intersection $\mathcal{O}=\bigcap_{n=1}^\infty \mathcal{O}_{x_n}$ is still a residual set and satisfies the desired property.   	
\end{proof}

%\,
%
%\noindent{\textbf{Declarations}}
%
%\,
%
%\noindent{\textbf{Supplementary data and Conflict of interest statements}} This document does not have any associated supplementary data. The authors have no Conflict of interest to declare.

%%%%%%%%%%%%%%%%%%%%%%%%%%%%%%%%%%%%%%%%%%%%%%%%%%%%%%%%%%%%%%%%%%
\appendix
\section{Weak KAM Theory and Aubry–Mather Theory}\label{sect_weakKAMandAubrymather}
This section provides some useful results of weak KAM theory and Aubry–Mather theory that are needed in this paper. Let $M$ be a closed connected smooth manifold. Throughout this section, ${\mathfrak H}:T^*M\to\R$ is assumed to be a continuous Hamiltonian that satisfies the following two conditions: 
\begin{itemize}
\item[(${\mathfrak H}1$)] (convexity) For every given $x\in M$, ${\mathfrak H}(x,p)$ 
is convex in $p\in T^*_xM$.
 \item[(${\mathfrak H}2$)]  (superlinearity) $\displaystyle{\frac{{\mathfrak H}(x,p)}{\|p\|_x}\to +\infty \quad \text{as~ $\|p\|_x\to +\infty$.} }$
\end{itemize}
The Hamilton-Jacobi equation associated to ${\mathfrak H}$ is as follows:
\begin{equation}\label{HJSansu}
{\mathfrak H}(x,d_xu)=c \quad\text{in $M$.}
\end{equation} 
Due to condition (${\mathfrak H}2$), it is well known (see e.g \cite{Barles_book}) that if $u$ is a viscosity subsolution of equation \eqref{HJSansu}, then $u$ is Lipschitz continuous with  Lipschitz constant $K_c$ given by
$$K_c=\sup\{\|p\|_x ~|~{\mathfrak H}(x,p)\leq c\}.$$

The {\em critical value} $c({\mathfrak H})$, also referred to as {\em ergodic constant} in the PDE literature, of the Hamiltonian ${\mathfrak H}$ is given by
\[c({\mathfrak H})=\inf_{f\in C^{1}(M, \R)}\max_{x\in M} ~{\mathfrak H}(x, d_xf).\]
%\[c({\mathfrak H})=\inf\{c\in \R : \text{equation \eqref{HJSansu} admits viscosity subsolutions}\}.\]
It was shown in \cite{Lions_Papanicolaou_Varadhan1987} (see also \cite{Fathi2012}) that $c=c({\mathfrak H})$ is the only constant for which equation \eqref{HJSansu} admits viscosity solutions.

The Lagrangian   ${\mathfrak L}:TM\to\R$ associated to the Hamiltonian ${\mathfrak H}$ is 
$${\mathfrak L}(x,v):=\sup_{p\in T^*_xM} p(v) -
{\mathfrak H}(x,p).
$$
${\mathfrak L}$ enjoys properties similar to ${\mathfrak H}$, namely, it is  continuous and satisfies:
\begin{itemize}
\item[(${\mathfrak L}1$)] (convexity) For every given $x\in M$,  
${\mathfrak L}(x,v)$ is convex in $v\in T_xM$.
 \item[(${\mathfrak L}2$)] (superlinearity) $\displaystyle{\frac{{\mathfrak L}(x,v)}{\|v\|_x}\to +\infty \quad \text{as~ $\|v\|_x\to +\infty$.} }$
\end{itemize}

For $t>0$, we define the minimal action function $h_t:M\times M\to \R$ by
\begin{equation*}\label{def h_t}
h_t(x,y)=\inf_{\gamma\in \Gamma^t_{x,y}}\int_0^t  \mathfrak{L}(\gamma(s),\dot\gamma(s))\, ds 
\end{equation*}
where $\Gamma^t_{x,y}$ is the family of absolutely continuous curves $\gamma:[0,t]\to M$ satisfying $\gamma(0)=x$ and $\gamma(t)=y$.

\begin{defn}[Peierls barrier]\label{dsadw14}
The {\em Peierls barrier} of $\mathfrak{L}$ is the function $h:M\times M\to\R$ defined by
\begin{equation*}
h(x,y):=\liminf_{t\to +\infty} \big[h_t(x,y)+c(\mathfrak{H})t\big].
\end{equation*}
\end{defn}

It satisfies the following properties, see \cite{Fathi_Siconolfi2005, Davini_Zavidovique2013}:
\begin{Pro}[Properties of $h$]\label{prop h}\ \\
{\rm (i)} The Peierls barrier
 $h$ is finite valued and Lipschitz
continuous.\\
\noindent{\rm (ii)} If $w$ is a viscosity subsolution of
${\mathfrak H}(x,d_xu)=c({\mathfrak H})$, then 
$$w(x)-w(y)\leq h(y,x)\quad\text{ for every $x,y\in M$}.$$
{\rm (iii)} For every fixed $y\in
    M$, the function $h(y,\cdot)$ is a viscosity solution of
${\mathfrak H}(x,d_xu)=c({\mathfrak H})$, and the function  $-h(\cdot,y)$ is a viscosity subsolution of
${\mathfrak H}(x,d_xu)=c({\mathfrak H})$.
\end{Pro}

\begin{defn}[Projected Aubry set]
Denote
$$  \mathcal{A}_{\mathfrak{L}}=\{x\in M \mid h(x,x)=0\},$$
we say that $\mathcal{A}_{\mathfrak{L}}$ is the {\em projected Aubry set} for the Lagrangian ${\mathfrak L}$. 
\end{defn}

The following result indicates that the projected Aubry set is a uniqueness set for the H-J equation, see for instance \cite{Fathi2012}.
\begin{Pro}\label{Aubry_uniqueness}
Let $w_-$ and $w_+$ be respectively viscosity subsolution and supersolution of the equation $\mathfrak{H}(x, d_xu)=c(\mathfrak{H})$. If $w_-(x)\leq w_+(x)$ for all $x\in\mathcal{A}_{\mathfrak{L}}$, then $w_-\leq w_+$  on $M$.	
\end{Pro}

We now proceed to introduce the Mather (minimizing) measures for the Lagrangian $\mathfrak L$. Mather measures, named after J. Mather, are closed probability measures that minimize the action of the Lagrangian.
Some good references are \cite{Mather1991, Mane1996, Fathi_Siconolfi2005}. 
  
\begin{defn}[Closed measure]\label{def closed measure}
A Borel probability measure $\tilde\mu$ on $TM$ is {\em closed} if 
\begin{itemize}
\item[\rm (1)]  $\displaystyle{\int_{TM} \|v\|_x\,d \tilde\mu(x,v)<\infty}$;
\item[\rm (2)] for every  $\varphi\in C^1(M)$, one has $\displaystyle{\int_{TM} d_x\varphi(v)\, d \tilde\mu (x,v) =0}$. 
\end{itemize}
\end{defn}

The relation linking closed probability measures to the critical value $c(\mathfrak H)
$ is as follows:
\begin{Pro}\label{Propertiesclosedmeas}
The following holds:
\begin{equation*}
	\min_{\tilde\mu}\int_{TM}{\mathfrak L}(x,v)\, d\tilde\mu(x,v)= -c(\mathfrak H)
\end{equation*}
where $\tilde\mu$ varies in the set of closed probability measures on $TM$. Hence, 
	\[
		\int_{TM}{\mathfrak L}(x,v)+c({\mathfrak H})\, d\tilde\mu(x,v) \geq 0,
	\]
	for any closed measure $\tilde\mu$ on $TM$.
\end{Pro}

\begin{defn}[Mather measure]\label{Def_Mather}
A {\em Mather measure} for the Lagrangian ${\mathfrak L}$ is a closed probability measure $\tilde\mu$ on $TM$ such that 
$$\int_{TM}{\mathfrak L}(x,v)\, d\tilde\mu(x,v)=-c({\mathfrak H}).$$
 \end{defn}
 
We remark that all Mather measures have support in a common compact subset of $TM$.  Denote by $\wtM({\mathfrak L})$ the family of all Mather measures on $TM$ for ${\mathfrak L}$. This set is convex, and compact in the weak-* topology.

\begin{defn}[Mather set]\label{Def_Matherset}
	The {\em Mather set} for the Lagrangian $\mathfrak L$ is defined as
\[\widetilde{\mathcal M}_{\mathfrak L}:=\overline{\bigcup_{\tilde\mu\in \wtM({\mathfrak L})}\supp \tilde\mu},\]
where $\supp \tilde\mu$ is the support of the measure $\tilde\mu$.
 The {\em projected Mather set} is defined by ${\mathcal M}_{\mathfrak L}=\pi \circ \widetilde{\mathcal M}_{\mathfrak L}$, where $\pi:TM\to M$ is the canonical projection.
\end{defn}
 The  Mather set  $\widetilde{\mathcal M}_{\mathfrak L}$ is a nonempty, compact subset of $TM$. Meanwhile, it is a classical result in Aubry-Mather theory that $ {\mathcal M}_{\mathfrak L}\subset {\mathcal A}_{\mathfrak L}$, and $\mathcal A_{\mathfrak L}\backslash{\mathcal M}_{\mathfrak L}\neq \emptyset$ could happen in general. In the case of a $C^2$ Tonelli Lagrangian, the Mather measures are invariant under the Euler-Lagrange flow, and the Mather set is also a flow invariant set. In \cite{Mather1991} Mather proved  the celebrated graph theorem:

\begin{The}\label{appendthm01}
	Let $\mathfrak{L}: TM\to \R$ be a $C^2$ Tonelli Lagrangian. Then, the Mather set $\widetilde{\mathcal M}_{\mathfrak L}$ is 	a Lipschitz graph above the projected Mather set ${\mathcal M}_{\mathfrak L}$, i.e., the restriction of the projection $\pi: TM\to M$ to 
	$\widetilde{\mathcal M}_{\mathfrak L}$ is injective with Lipschitz inverse.\end{The}

For a measure $\tilde\nu$ on $TM$, we denote by $\nu$ its projection $\pi_\#\,\tilde\nu$ on $M$ where $\pi:TM\to M$ is the canonical projection. Namely, for each Borel subset $A\subset M$,
\[\nu(A):=\tilde\nu(\pi^{-1}(A)).\]
In other words, 
\[ \int_M f(x)\,d\mu(x)=\int_{TM} (f\circ\pi)(x,v)\,d\tilde\mu(x,v),\quad\text{for every $f\in C(M)$}.\]

\begin{defn}[Projected Mather measure]
A {\em projected Mather measure} is a probability measure on $M$ of the form $\mu=\pi_\#\tilde\mu$, where $\tilde\mu$ is a Mather measure for ${\mathfrak L}$.
\end{defn}

\begin{Rem}
Let  $\mathfrak{M}({\mathfrak L})$ be the family of all projected Mather measures on $M$ for ${\mathfrak L}$, then the graph Theorem  \ref{appendthm01} implies that in the case of a $C^2$ Tonelli Lagrangian, the linear map $\pi_\#: \wtM({\mathfrak L})\to \mathfrak{M}({\mathfrak L}) $ is injective. 
\end{Rem}

%\noindent\textbf{Acknowledgments.} 

%\bibliography{mybib}
%\bibliographystyle{plain}%{alpha}

\end{document}